\documentclass[12pt,letterpaper,reqno]{amsart}

\usepackage{mathptmx}
\usepackage{amsmath}
\usepackage{amssymb}
\usepackage{mathtools}
\usepackage{mathrsfs}
\usepackage{bm}
\usepackage{enumitem}
\usepackage{microtype}
\usepackage{xcolor}

\usepackage[capposition=top]{floatrow}
\usepackage{euscript,eufrak,verbatim}
\usepackage{graphicx}
\usepackage{appendix}
\usepackage{amsthm}
\usepackage[all]{xy}
\usepackage{bbm}
\usepackage{amscd}

\usepackage[
  colorlinks=true,
  linkcolor=red,
  citecolor=blue,
  urlcolor=blue,
  pdfencoding=auto,
  psdextra
]{hyperref}
\usepackage[nameinlink,noabbrev,capitalise]{cleveref}

\theoremstyle{plain}
\newtheorem*{main theorem}{Main Theorem}
\newtheorem{theorem}{Theorem}[section]
\newtheorem{lemma}[theorem]{Lemma}
\newtheorem{proposition}[theorem]{Proposition}
\newtheorem{corollary}[theorem]{Corollary}

\newtheorem{question}[theorem]{Question}
\newtheorem{claim}[theorem]{Claim}

\theoremstyle{definition}
\newtheorem{definition}[theorem]{Definition}

\theoremstyle{remark}
\newtheorem{remark}[theorem]{Remark}

\numberwithin{equation}{section}
\setlist[enumerate,1]{label=\textup{(\roman*)},leftmargin=*,itemsep=2pt}
\setlist[itemize]{leftmargin=*,itemsep=2pt}

\title[Lowering topological entropy for amenable group actions]{Lowering topological entropy and asymptotic $h$-expansiveness for amenable group actions}
\author[X. Wang]{Xiaochen Wang}
\address{School of Mathematical Sciences, University of Science and Technology of China,
Hefei, 230026, Anhui, People’s Republic of China}
\email{xcwang97@ustc.edu.cn}
\subjclass[2020]{Primary 37A35, 37B40; Secondary 37A15, 37B05}
\keywords{amenable group action, lowering topological entropy, Bowen topological entropy, D-hereditarily lowerable, hereditarily uniformly lowerable, asymptotically $h$-expansive, tail entropy}
\date{}

\hypersetup{
  pdftitle={Lowering topological entropy and asymptotic h-expansiveness for amenable group actions},
  pdfauthor={Xiaochen Wang},
  pdfsubject={Mathematics},
  pdfkeywords={amenable group action, lowering topological entropy, Bowen topological entropy, D-hereditarily lowerable, hereditarily uniformly lowerable, asymptotically h-expansive, tail entropy}
}

\begin{document}

\begin{abstract}
Let $G$ be a countably infinite discrete amenable group acting continuously on a compact metric space $X$.  
We study the problem of lowering topological entropy over subsets of $X$ and its connections with asymptotic $h$-expansiveness.
Let $\{F_n\}$ be a tempered F{\o}lner sequence such that $e_G\in F_1\subseteq F_2\subseteq\cdots$ and $|F_n|/\log n\to\infty$. 
If $(X,G)$ has finite topological entropy, then, for every $h\in[0,h_{\mathrm{top}}(X,G)]$,there exists a non-empty compact set $K_h\subseteq X$ whose topological entropy, Bowen topological entropy and packing topological entropy along $\{F_n\}$ are all equal to $h$.
  
We next consider D-hereditary lowerability, which requires every non-empty analytic (Souslin) subset to be lowerable with respect to Bowen topological entropy. 
We prove that every subshift over a finite alphabet is D-hereditarily lowerable along an increasing F{\o}lner sequence. 
Combining principal quasi-symbolic extensions with a dimensional entropy inequality for factor maps, we further prove that every asymptotically $h$-expansive $G$-system is D-hereditarily lowerable along an increasing F{\o}lner sequence.
More generally, suppose that $(X,G)$ has finite topological entropy and tail entropy $a$, and let $\{F_n\}$ be an increasing F{\o}lner sequence. If an analytic set $K\subseteq X$ has Bowen topological entropy $H>a$ along $\{F_n\}$, then, for every $h\in[0,H-a]$, there exists a non-empty set $K_h\subseteq K$ such that $h\le h_{\mathrm{top}}^B(K_h,{F_n})\le h+a$. 
Finally, for tempered increasing F{\o}lner sequences satisfying the above growth condition, we prove that asymptotic $h$-expansiveness is equivalent to hereditary uniform lowerability.  
\end{abstract}

\maketitle

\pagestyle{headings}
\markboth{}{}

\section{Introduction}\label{section: introduction}
Throughout the paper, by a $G$-system $(X,G)$ we mean a compact metric space $X$ and a countably infinite discrete amenable group $G$ acting continuously on $X$ by homeomorphisms.
A basic problem in entropy theory is whether a prescribed amount of entropy below that of a given system can be realized by passing to an appropriate factor or subset. 
Its measure-theoretic prototype is Sinai's factor theorem, i.e., every ergodic measure-preserving $\mathbb Z$-system of positive entropy has Bernoulli factors of every smaller entropy \cite{Sinai1962}. 
In the topological category, Shub and Weiss \cite{ShubWeiss1991} asked whether every positive-entropy system admits a continuous factor whose entropy lies strictly between zero and that of the original system.  Lindenstrauss \cite{Lindenstrauss1995} proved that, for finite-dimensional systems, every value between zero and the system's entropy is realized by a factor.
Later, analogous results for extensions of non-trivial minimal systems were established using mean dimension and the small boundary property \cite{Lindenstrauss1999}. Motivated by this circle of questions, Huang, Ye, and Zhang \cite{Huang2010YeZhang,Huang2014YeZhang} replaced factors by subsets and initiated the systematic study of entropy lowering for $\mathbb Z$-systems.  They established, among other results, lowerability for finite-entropy systems, characterized asymptotic $h$-expansiveness in terms of hereditary uniform lowerability, and proved that asymptotically $h$-expansive systems are D-hereditarily lowerable.  Unlike a factor, a subset need not be invariant and hence need not carry an induced dynamical system.  Entropy lowering over subsets therefore concerns the distribution of orbit complexity within the phase space, rather than the construction of quotient systems.

We consider this problem for a $G$-system. 
For $\mathbb Z$-actions, the topological entropy of compact subsets is defined using fixed-length orbit covers, while Bowen \cite{Bowen1973} introduced an entropy for arbitrary subsets through a Carath\'eodory construction modelled on Hausdorff dimension. 
Feng and Huang \cite{FengHuang2012} subsequently introduced packing topological entropy as the dynamical counterpart of packing dimension.  
In the setting of amenable group actions, orbit intervals are replaced by finite sets taken from a F{\o}lner sequence.  
Zheng and Chen \cite{ZhengChen2016} introduced Bowen topological entropy along a fixed F{\o}lner sequence, while Dou, Zheng, and Zhou \cite{Dou2023ZhengZhou} introduced the corresponding packing version.  
The Ornstein–Weiss lemma \cite{OrnsteinWeiss1987} ensures that, for the whole space, the usual topological entropy is independent of the chosen Følner sequence.
Under appropriate temperedness and growth assumptions, both the Bowen entropy and the packing entropy of the whole space agree with this quantity \cite{ZhengChen2016,Dou2023ZhengZhou}. 
This independence need not persist for non-invariant subsets, since the translation invariance required by the Ornstein--Weiss argument is then lost.  
Indeed, the symbolic construction in \cite[Proposition~5.8]{Dou2023ZhengZhou} expresses the Bowen and packing topological entropies of a compact subset of the full shift in terms of the lower and upper densities, respectively, of a prescribed coordinate set along the chosen F{\o}lner sequence.  
Changing the sequence can therefore change both entropy values.  
Thus entropy lowering over subsets must be formulated relative to a fixed averaging sequence.

Against this background, we investigate which forms of entropy lowering survive for amenable group actions and how they are governed by the entropy structure of the system.  
This is not a formal replacement of the intervals $\{0,\ldots,n-1\}$ by F{\o}lner sets.  
The arguments for a single transformation can exploit the order and nested geometry of orbit intervals, whereas a general F{\o}lner sequence has no canonical order and may consist of sets with substantially different geometry.  
Moreover, the fixed-length covering procedure, the Carath\'eodory construction, and the packing construction need not yield the same entropy on a non-invariant subset. 
We therefore ask whether, along a F{\o}lner sequence, every prescribed entropy value can be realized by a compact subset simultaneously for all three notions; whether prescribed Bowen entropy value can be realized inside a given Souslin set; and how hereditary forms of lowerability are related to asymptotic $h$-expansiveness and tail entropy.  
These questions identify the roles played by symbolic models and tail complexity in entropy lowerability and clarify which features of the $\mathbb Z$-theory are genuinely independent of the ordered structure of time.

The results of this paper are organized in three parts.  
First, along a tempered F{\o}lner sequence satisfying the standard superlogarithmic growth condition, we prove that every finite-entropy $G$-system admits, for each prescribed entropy level, a compact subset on which the topological, Bowen, and packing topological entropies coincide with that level. We also establish hereditary uniform lowerability for expansive actions along an appropriate regular class of F{\o}lner sequences. Second, we turn to hereditary lowering for Bowen topological entropy. We first prove that finite-alphabet subshifts are D-hereditarily lowerable along every increasing F{\o}lner sequence and then, using principal quasi-symbolic extensions, extend this conclusion to asymptotically $h$-expansive amenable actions. More generally, if a finite-entropy system has tail entropy $a$ and a Souslin set $K$ has Bowen entropy $H>a$, then every $0\leq h\leq H-a$ can be realized inside $K$ up to an error of at most $a$. Finally, for a regular class of F{\o}lner sequences, we prove that asymptotic $h$-expansiveness is equivalent to hereditary uniform lowerability.  
Taken together, these results identify symbolic models as the source of exact hereditary lowering in the zero-tail entropy case and show that, in general, tail entropy quantifies the obstruction to exact entropy lowerability.

\subsection*{Main results}
Our first result shows that every prescribed entropy value of a finite-entropy system can be realized on a single compact subset simultaneously for all three notions of subset entropy.

\begin{theorem}\label{thm-main-finite}
Let $(X,G)$ be a $G$-system with finite topological entropy, and let
$\{F_n\}$ be a tempered F{\o}lner sequence such that
\[
 e_G\in F_1\subseteq F_2\subseteq\cdots,
 \qquad
 \lim_{n\to\infty}\frac{|F_n|}{\log n}=+\infty.
\]
Then for every $0\le h\le h_{\mathrm{top}}(X,G)$ there exists a nonempty
compact set $K_h\subseteq X$ such that
\[
 h_{\mathrm{top}}^B(K_h,\{F_n\})
 =h_{\mathrm{top}}^P(K_h,\{F_n\})
 =h_{\mathrm{top}}(K_h,\{F_n\})
 =h.
\]
In particular, $(X,G)$ is lowerable, D-lowerable, and P-lowerable along
$\{F_n\}$.
\end{theorem}

Requiring prescribed entropy value to be realized inside every Souslin subset leads to a stronger form of lowering. 
The relevant structural quantity is topological tail entropy. Recall that $(X,G)$ is asymptotically
$h$-expansive if its topological tail entropy $h_{\mathrm{top}}^*(X,G)=0$. Our second main theorem shows that this zero-tail condition is sufficient for D-hereditary lowerability.

\begin{theorem}\label{thm-main-DHL}
Let $(X,G)$ be a $G$-system with finite topological entropy, and let
$\{F_n\}$ be a F{\o}lner sequence such that
\[
 e_G\in F_1\subseteq F_2\subseteq\cdots.
\]
If $(X,G)$ is asymptotically $h$-expansive, then it is D-hereditarily
lowerable along $\{F_n\}$.
\end{theorem}

The preceding theorem treats the zero-tail entropy case. More generally,
topological tail entropy controls the error with which a prescribed Bowen
entropy can be realized inside a Souslin subset.

\begin{theorem}\label{thm-main-tailinterval}
Let $(X,G)$ be a $G$-system with finite topological entropy, let $\{F_n\}$ be a
F{\o}lner sequence satisfying
$e_G\in F_1\subseteq F_2\subseteq\cdots$, and let $K\subseteq X$ be a nonempty
Souslin set.  Put
\[
 H=h_{\mathrm{top}}^B(K,\{F_n\}),
 \qquad
 a=h_{\mathrm{top}}^*(X,G).
\]
If $H>a$, then for every $0\le h\le H-a$ there exists
$K_h\subseteq K$ such that
\[
 h\le h_{\mathrm{top}}^B(K_h,\{F_n\})\le h+a.
\]
\end{theorem}

Thus tail entropy gives an explicit error bound for lowering Bowen entropy
inside a Souslin set.  When $a=0$, the two bounds coincide and one recovers
D-hereditary lowerability.  In this sense, topological tail entropy
quantitatively controls the obstruction to exact hereditary lowering for
Bowen entropy.

For compact subsets and topological entropy, we obtain the stronger property of hereditary uniform lowerability.  
Along the regular F{\o}lner
sequences used in \eqref{FS-condition} below, this property characterizes
asymptotic $h$-expansiveness.

\begin{theorem}\label{thm-main-HUL}
Let $(X,G)$ be a $G$-system and let $\{F_n\}$ be a tempered
F{\o}lner sequence satisfying \eqref{FS-condition}.  
Then $(X,G)$ is
asymptotically $h$-expansive if and only if it is hereditarily uniformly
lowerable along $\{F_n\}$.
\end{theorem}

Theorems~\ref{thm-main-DHL} and~\ref{thm-main-HUL} describe two complementary
hereditary lowering properties.  
The D-hereditary statement concerns Bowen
topological entropy on arbitrary Souslin sets and requires only an increasing
F{\o}lner sequence, whereas hereditary uniform lowerability concerns the
topological entropy of compact subsets and produces countable compact subsets with at most one limit point. 
Both are governed by the disappearance of tail complexity.

\subsection*{Organization of the paper}

The paper is organized as follows.  Section~\ref{preliminary} introduces the three notions of entropy for subsets and the associated lowering properties.  We also recall the definitions and results concerning topological tail entropy and extensions that will be used later.  Section~\ref{section-Finite-entropy-system} proves the simultaneous lowering theorem for finite-entropy systems, namely Theorem~\ref{thm-main-finite}, using distribution principles and relative entropy estimates.  The same section also treats expansive actions and establishes their hereditary uniform lowerability under the stated assumptions.  Section~\ref{section-DHL-ahe} studies the relation between Bowen topological entropy and Hausdorff dimension.  We first prove D-hereditary lowerability for finite-alphabet subshifts and then use principal quasi-symbolic extensions to prove Theorem~\ref{thm-main-DHL}.  The section concludes with the quantitative estimate involving topological tail entropy given in Theorem~\ref{thm-main-tailinterval}.  Finally, Section~\ref{section-asyhexpansive=hul} proves the equivalence between asymptotic $h$-expansiveness and hereditary uniform lowerability in Theorem~\ref{thm-main-HUL}, and establishes the relevant preservation results under principal extensions.

\section{Preliminaries}\label{preliminary}
Let $(X,G)$ be a $G$-system, where $X$ is a compact metric space equipped with metric $d$, and $G$ is a countably infinite discrete amenable group acting continuously on $X$ via homeomorphisms. 
Denote by $\mathcal{F}(G)$ the set of all nonempty finite subsets of $G$. Recall that a countable discrete group $G$ is \emph{amenable} if it admits a \emph{F{\o}lner sequence} $\{F_n\}$, i.e., each $F_n$ is a nonempty finite subset of $G$ and for every $s \in G$, 
\[
\frac{|F_n \bigtriangleup s F_n|}{|F_n|} \to 0 \quad \text{as } n \to \infty.
\]
Moreover, a sequence $\{F_n\}$ is called \emph{tempered} if there exists $C > 0$ such that 
\[
\left|\bigcup_{k < n} F_k^{-1} F_n \right| \leq C |F_n| \text{ for all } n > 1.    
\] 

Let $(X,G)$ be a $G$-system, $K \subseteq X$, and $\mathcal{U}$ a family of subsets of $X$. Denote by $\mathrm{diam}(K)$ the diameter of $K$, and set $\|\mathcal{U}\| = \sup \{ \mathrm{diam}(U) : U \in \mathcal{U} \}$. 
We write $K \succeq \mathcal{U}$ if $K \subseteq U$ for some $U \in \mathcal{U}$; otherwise, $K \nsucceq \mathcal{U}$. 
If $\mathcal{W}$ is another family of subsets of $X$, we say $\mathcal{W}$ is \emph{finer} than $\mathcal{U}$ (denoted $\mathcal{W} \succeq \mathcal{U}$) if $W \succeq \mathcal{U}$ for every $W \in \mathcal{W}$.

A \emph{cover} of $X$ is a finite family of Borel subsets of $X$ whose union is $X$, and a \emph{partition} of $X$ is a cover whose elements are pairwise disjoint. 
Denote by $\mathcal{C}_X$ (respectively, $\mathcal{C}_X^o$, $\mathcal{P}_X$) the set of all finite covers (respectively, open covers, measurable partitions) of $X$.

Note that if $\mathcal{U} \in \mathcal{C}_X^o$, then it has a \textit{Lebesgue number} $\lambda > 0$, and $\mathcal{W} \succeq \mathcal{U}$ whenever $\|\mathcal{W}\| < \lambda$. 
For $\alpha \in \mathcal{P}_X$ and $x \in X$, let $\alpha(x)$ denote the unique element of $\alpha$ containing $x$. 
Given $\mathcal{U}_1, \mathcal{U}_2 \in \mathcal{C}_X$, define
\[
\mathcal{U}_1 \vee \mathcal{U}_2 = \{ U_1 \cap U_2 : U_1 \in \mathcal{U}_1,~U_2 \in \mathcal{U}_2\}.  
\] 
Then $\mathcal{U}_1 \vee \mathcal{U}_2 \in \mathcal{C}_X $ and $\mathcal{U}_1 \vee \mathcal{U}_2 \succeq \mathcal{U}_1 $. Moreover, $ \mathcal{U}_1 \succeq \mathcal{U}_2 $ if and only if $ \mathcal{U}_1 \vee \mathcal{U}_2 \succeq \mathcal{U}_1 $ and $ \mathcal{U}_1 \succeq \mathcal{U}_1 \vee \mathcal{U}_2$. 
For each $\mathcal{U} \in \mathcal{C}_X$ and finite subset $F \subseteq G$, define
\[
\mathcal{U}_F = \bigvee_{g \in F} g^{-1}\mathcal{U}.   
\]

\subsection{Topological Entropy for G-Systems}
Let $(X,G)$ be a $G$-system, $K \subseteq X$, and $\mathcal{U} \in \mathcal{C}_X^o$. Define $N(\mathcal{U}, K)$ as the minimal cardinality of subfamilies $\mathcal{V} \subseteq \mathcal{U}$ such that $\bigcup \mathcal{V} \supseteq K$, where $\bigcup \mathcal{V} = \bigcup_{V \in \mathcal{V}} V$. By convention, set $N(\mathcal{U}, \emptyset) = 0$ and $N(\mathcal{U}) = N(\mathcal{U}, X)$.

For every finite open cover $\mathcal{U}$ of $X$, the function $F \mapsto \log N(\mathcal{U}_F)$ on the collection of finite subsets of $G$ is right $G$-invariant and subadditive, i.e., 
\[
\log N(\mathcal{U}_{E \cup F}) \le \log N(\mathcal{U}_E) + \log N(\mathcal{U}_F).
\]
Consequently, the quantity 
\[
\frac{1}{|F|} \log N(\mathcal{U}_F)
\]
converges to a limit as $F$ becomes increasingly invariant. This limit is independent of the choice of F{\o}lner sequence $\{F_n\} $ for $G$. The \emph{topological entropy of $\mathcal{U}$} is defined as
\[
h_{\mathrm{top}}(\mathcal{U}) = \lim_{n \to \infty} \frac{1}{|F_n|} \log N(\mathcal{U}_{F_n}),
\]
and the \emph{topological entropy of $(X,G)$} is
\[
h_{\mathrm{top}}(X, G) = \sup_{\mathcal{U}} h_{\mathrm{top}}(\mathcal{U}),
\]
where the supremum is taken over all finite open covers $\mathcal{U}$ of $X$.
Note that for any $s \in G$ and any finite subset $F \subseteq G$,
\[
N(\mathcal{U}_{Fs}, K) = N(s^{-1}\mathcal{U}_F, K) = N(\mathcal{U}_F, sK).
\]
Now define
\[
h_{\mathrm{top}}(K, \mathcal{U}, \{F_n\}) = \limsup_{n \to \infty} \frac{1}{|F_n|} \log N(\mathcal{U}_{F_n}, K).
\]
Clearly, $h_{\mathrm{top}}(K, \mathcal{U}, \{F_n\})$ is increasing in $\mathcal{U}$: if $\mathcal{V} \succeq \mathcal{U}$, then
\[
h_{\mathrm{top}}(K, \mathcal{U}, \{F_n\}) \le h_{\mathrm{top}}(K, \mathcal{V}, \{F_n\}).
\]
The \emph{topological entropy of $K$} is defined as
\[
h_{\mathrm{top}}(K, \{F_n\}) = \sup_{\mathcal{U} \in \mathcal{C}_X^o} h_{\mathrm{top}}(K, \mathcal{U}, \{F_n\}).
\]
	
Let $(X,G)$ be a $G$-system and define the \emph{hyperspace} of $X$ by
\[
2^X = \{ K \subseteq X : K \text{ is a non-empty compact subset} \}.
\]
This space is endowed with the Hausdorff metric $d_H$ given by
\begin{align*}
d_H(K,L) &= \max\left\{ \max_{x \in K} d(x,L),\ \max_{x \in L} d(x,K) \right\} \\
&= \inf\left\{ \varepsilon > 0 : K \subseteq U_{\varepsilon}(L) \text{ and } L \subseteq U_{\varepsilon}(K) \right\},
\end{align*}
where $U_{\varepsilon}(K) = \{ x \in X : d(x,K) < \varepsilon \}$. 

The \emph{entropy hyper-function} $H : 2^X \to [0, h_{\mathrm{top}}(X,G)]$ of $(X,G)$ is defined by
\[
H(K) = h_{\mathrm{top}}(K, \{F_n\}), \quad K \in 2^X.
\]

We then have the following properties.

\begin{proposition}\label{prop-coverentropy-property}
Let $(X,G)$ be a $G$-system, $\{F_n\} $ a F{\o}lner sequence of $G$, and $\mathcal{U} \in \mathcal{C}_X^o$. Then:
\begin{enumerate}
\item[(1)] For any $K \subseteq X$ and $s \in G$,
\[
h_{\mathrm{top}}(sK, \mathcal{U}, \{F_n\}) = h_{\mathrm{top}}(K, \mathcal{U}, \{F_n s\}).
\]
\item[(2)] The function 
\[
h_{\mathrm{top}}(\cdot, \mathcal{U}, \{F_n\}) : 2^X \to [0, \infty)
\]
is Borel measurable.
\item[(3)] The entropy hyper-function $H(\cdot) = h_{\mathrm{top}}(\cdot, \{F_n\})$ is Borel measurable.
\end{enumerate}
\end{proposition}

\subsection{Sep, Spn-Entropy for $G$-Systems}
Let $(X,G)$ be a $G$-system, $d$ be the metric on $X$, $\{F_n\} $ be a F{\o}lner sequence of $G$ and $K \subseteq X$. 
For each nonempty finite subset $F$ of $G$, define the metric $d_F$ on $X$ by 
\[
d_F(x,y) = \max_{s \in F} d(sx,sy).
\]
Let $F$ be a nonempty finite subset of $G$ and $\varepsilon > 0$. A set $D \subseteq K$ is called an \emph{$(d,F,\varepsilon)$-separated set} of $K$ if $d_F(x,y) \geq \varepsilon$ for all distinct $x,y \in D$, and an \emph{$(d,F,\varepsilon)$-spanning set} of $K$ if for every $x \in K$ there exists $y \in D$ such that $d_F(x,y) < \varepsilon$. 
Denote by $\mathrm{sep}(d,F,\varepsilon,K)$ the maximum cardinality of an $(d,F,\varepsilon)$-separated set of $K$, and by $\mathrm{spn}(d,F,\varepsilon,K)$ the minimum cardinality of an $(d,F,\varepsilon)$-spanning set of $K$.
Define
\begin{align*}
h_{\mathrm{sep}}(d,\varepsilon,K) &= \limsup_{n\to \infty} \frac{1}{|F_n|} \log \mathrm{sep}(d,F_n,\varepsilon,K), \\
h_{\mathrm{spn}}(d,\varepsilon,K) &= \limsup_{n\to \infty} \frac{1}{|F_n|} \log \mathrm{spn}(d,F_n,\varepsilon,K),
\end{align*}
and
\begin{align*}
h_{\mathrm{sep}}(d,K) &= \lim_{\varepsilon \to 0^+} h_{\mathrm{sep}}(d,\varepsilon,K), \\
h_{\mathrm{spn}}(d,K) &= \lim_{\varepsilon \to 0^+} h_{\mathrm{spn}}(d,\varepsilon,K).
\end{align*}
For convenience, we denote $\mathrm{sep}(d,F,\varepsilon) = \mathrm{sep}(d,F,\varepsilon,X)$, $h_{\mathrm{sep}}(d,\varepsilon) = h_{\mathrm{sep}}(d,\varepsilon,X)$, and $h_{\mathrm{sep}}(d) = h_{\mathrm{sep}}(d,X)$ (respectively, $\mathrm{spn}(d,F,\varepsilon)$, $h_{\mathrm{spn}}(d,\varepsilon)$, and $h_{\mathrm{spn}}(d)$).

We have the following results.
\begin{proposition}    
Let $(X,G)$ be a $G$-system, $\{F_n\} $ a F{\o}lner sequence of $G$, $K \subseteq X$, and $\varepsilon > 0$. Then    
\begin{enumerate}		
    \item[(1)] $\mathrm{spn}(d,F_n,\varepsilon,K) \leq \mathrm{sep}(d,F_n,\varepsilon,K) \leq \mathrm{spn}(d,F_n,\varepsilon/2,K)$.		
    
    \item[(2)] $h_{\mathrm{spn}}(d,K) = h_{\mathrm{sep}}(d,K)$, and this quantity is independent of the choice of a compatible metric $d$ on $X$.		
    
    \item[(3)] If $\mathcal{U} \in \mathcal{C}_X^o$ has a Lebesgue number $\delta > 0$, then for any $\delta' \in (0, \delta/2)$ and every $\mathcal{V} \in \mathcal{C}_X^o$ with $\|\mathcal{V}\| \leq \delta'$, one has 
    \[
    N(\mathcal{U}_{F_n}, K) \leq \mathrm{spn}(d, F_n, \delta', K) \leq \mathrm{sep}(d, F_n, \delta', K) \leq N(\mathcal{V}_{F_n}, K).
    \]		
    
    \item[(4)] If $\{\mathcal{U}_n\}_{n \in \mathbb{N}}$ is a sequence of finite open covers of $X$ satisfying $\|\mathcal{U}_n\| \to 0$ as $n \to \infty$, then 
    \[
    h_{\mathrm{spn}}(d,K) = h_{\mathrm{sep}}(d,K) = \lim_{n \to \infty} h_{\mathrm{top}}(K, \mathcal{U}_n, \{F_n\}) = h_{\mathrm{top}}(K, \{F_n\}).
    \]    
\end{enumerate}
\end{proposition}

\subsection{Bowen Topological Entropy for $G$-Systems}
Let $(X,G)$ be a $G$-system, $\{F_n\} $ a F{\o}lner sequence of $G$, $\mathcal{U} \in \mathcal{C}_X^o$, and $K \subseteq X$. For each $n \geq 1$, define
\[
\mathcal{W}_{F_n}(\mathcal{U}) = \left\{ \mathbb{U} = \{U_g\}_{g \in F_n} : U_g \in \mathcal{U},\ \forall g \in F_n \right\}.
\]
For $\mathbb{U} \in \mathcal{W}_{F_n}(\mathcal{U})$, the integer $m(\mathbb{U}) := |F_n|$ is called the \emph{length} of $\mathbb{U}$, and we define
\[
X(\mathbb{U}) = \bigcap_{g \in F_n} g^{-1} U_g = \left\{ x \in X : gx \in U_g,\ \forall g \in F_n \right\}.
\]
A subset $\Lambda \subseteq \bigcup_{n \geq 1} \mathcal{W}_{F_n}(\mathcal{U})$ is said to \emph{cover} $K$ if
\[
K \subseteq \bigcup_{\mathbb{U} \in \Lambda} X(\mathbb{U}).
\]
	
For $s \ge 0$ and $N \in \mathbb{N}$, define
\[
\mathcal{M}(K, \mathcal{U}, N, s, \{F_n\}) = \inf \left\{ \sum_{\mathbb{U} \in \Lambda} e^{-s \cdot m(\mathbb{U})} : \Lambda \subseteq \bigcup_{j \geq N} \mathcal{W}_{F_j}(\mathcal{U}) \text{ covers } K \right\}.
\]
Note that $\mathcal{M}(K, \mathcal{U}, N, s, \{F_n\})$ increases with $N$, and define
\[
\mathcal{M}(K, \mathcal{U}, s, \{F_n\}) = \lim_{N \to \infty} \mathcal{M}(K, \mathcal{U}, N, s, \{F_n\}).
\]
The \emph{Bowen topological entropy of $K$ relative to $\mathcal{U}$ (along the F{\o}lner sequence $\{F_n\} $)} is defined by
\begin{align*}
h_{\mathrm{top}}^B(K, \mathcal{U}, \{F_n\}) &= \inf \{ s \in \mathbb{R} : \mathcal{M}(K, \mathcal{U}, s, \{F_n\}) = 0 \} \\
&= \sup \{ s \in \mathbb{R} : \mathcal{M}(K, \mathcal{U}, s, \{F_n\}) = \infty \}.
\end{align*}
The \emph{Bowen topological entropy of $K$ (along the F{\o}lner sequence $\{F_n\} $)} is defined by
\[
h_{\mathrm{top}}^B(K, \{F_n\}) = \sup_{\mathcal{U} \in \mathcal{C}_X^o} h_{\mathrm{top}}^B(K, \mathcal{U}, \{F_n\}).
\]

The Bowen topological entropy admits an alternative definition. For $K \subseteq X$, $s \ge 0$, $N \in \mathbb{N}$, and $\varepsilon > 0$, define
\begin{align*}    
\mathcal{M}(K, N, \varepsilon, s, \{F_n\}) = \inf \Bigg\{ \sum_i & e^{-s \cdot |F_{n_i}|} : \{B_{F_{n_i}}(x_i, \varepsilon)\}_i \text{ is a finite or countable family} \\     
& \text{such that } x_i \in X,\ n_i \ge N,\text{ and } \bigcup_i B_{F_{n_i}}(x_i, \varepsilon) \supseteq K \Bigg\}.
\end{align*}
Since $\mathcal{M}(K, N, \varepsilon, s, \{F_n\})$ is non-decreasing in $N$ and non-increasing in $\varepsilon$, we define
\begin{align*}    
\mathcal{M}(K, \varepsilon, s, \{F_n\}) &= \lim_{N \to \infty} \mathcal{M}(K, N, \varepsilon, s, \{F_n\}), \\    
\mathcal{M}(K, s, \{F_n\}) &= \lim_{\varepsilon \to 0} \mathcal{M}(K, \varepsilon, s, \{F_n\}).
\end{align*}
The Bowen topological entropy of a subset $K$ along $\{F_n\}_{n \in \mathbb{N}}$, denoted $h_{\mathrm{top}}^B(K, \{F_n\})$, can be equivalently defined as the critical value of the parameter $s$ at which $\mathcal{M}(K, s, \{F_n\})$ jumps from $\infty$ to $0$, i.e.,
\[ 
\mathcal{M}(K, s, \{F_n\}) = 
\begin{cases}    
\infty, & s < h_{\mathrm{top}}^B(K, \{F_n\}), \\    
0, & s > h_{\mathrm{top}}^B(K, \{F_n\}).
\end{cases} 
\]
This equivalence is established in \cite[Page 74]{Pesin1997}.

We have the following results relating topological entropy and Bowen topological entropy \cite{ZhengChen2016, Dou2018Zhang}.

\begin{proposition}\label{prop-topoent-Bowentopoent}
Let $(X,G)$ be a $G$-system, $\{F_n\} $ a tempered F{\o}lner sequence of $G$ satisfying the growth condition $\lim\limits_{n \to \infty} \frac{|F_n|}{\log n} = \infty$, and let $K_1, K_2, \dots, K \subseteq X$ and $\mathcal{U} \in \mathcal{C}_X^o$. Then:
\begin{enumerate}	
    \item[(1)] $h_{\mathrm{top}}(X, G) = h_{\mathrm{top}}^B(X, \{F_n\})$.	
    
    \item[(2)] $h_{\mathrm{top}}^B\left(\bigcup_{n \geq 1} K_n, \mathcal{U}, \{F_n\}\right) = \sup_{n \geq 1} h_{\mathrm{top}}^B(K_n, \mathcal{U}, \{F_n\})$.
\end{enumerate}
\end{proposition}

\subsection{Packing Topological Entropy for $G$-Systems}
Let $(X,G)$ be a $G$-system and $\{F_n\} $ a F{\o}lner sequence of $G$. For $K \subseteq X$, $s \ge 0$, $N \in \mathbb{N}$, and $\varepsilon > 0$, define 
\begin{align*}    
P(K, N, \varepsilon, s, \{F_n\}) = \sup \Bigg\{ \sum_i & e^{-s \cdot |F_{n_i}|} : \{\overline{B}_{F_{n_i}}(x_i, \varepsilon)\}_i \text{ is a finite or countable pairwise}  \\     
& \text{disjoint family such that } x_i \in K \text{ and } n_i \ge N \Bigg\}.
\end{align*}
Since $P(K, N, \varepsilon, s, \{F_n\})$ is non-increasing in $N$, we define
\[ 
P(K, \varepsilon, s, \{F_n\}) = \lim_{N \to \infty} P(K, N, \varepsilon, s, \{F_n\}).
\]
Define
\[ 
\mathcal{P}(K, \varepsilon, s, \{F_n\}) = \inf \left\{ \sum_{i \ge 1} P(K_i, \varepsilon, s, \{F_n\}) : \bigcup_{i \ge 1} K_i \supseteq K \right\}.
\]
There exists a critical value of the parameter $s$, denoted by $h_{\mathrm{top}}^P(K, \varepsilon, \{F_n\})$, at which $\mathcal{P}(K, \varepsilon, s, \{F_n\})$ jumps from $\infty$ to $0$, i.e.,
\[ 
\mathcal{P}(K, \varepsilon, s, \{F_n\}) = 
\begin{cases}    
\infty, & s < h_{\mathrm{top}}^P(K, \varepsilon, \{F_n\}), \\    
0, & s > h_{\mathrm{top}}^P(K, \varepsilon, \{F_n\}).
\end{cases}
\]
Since $h_{\mathrm{top}}^P(K, \varepsilon, \{F_n\})$ is non-decreasing as $\varepsilon$ decreases, the \emph{packing topological entropy of $K$ (along the F{\o}lner sequence $\{F_n\} $)} is defined as
\[ 
h_{\mathrm{top}}^P(K, \{F_n\}) = \lim_{\varepsilon \to 0} h_{\mathrm{top}}^P(K, \varepsilon, \{F_n\}).
\]

We have the following results relating topological entropy and packing topological entropy \cite{Dou2023ZhengZhou}.

\begin{proposition}    
Let $(X,G)$ be a $G$-system and $\{F_n\} $ a F{\o}lner sequence of $G$ satisfying the growth condition $\lim\limits_{n \to \infty} \frac{|F_n|}{\log n} = \infty$. Then     
\[ 
h_{\mathrm{top}}(X, G) = h_{\mathrm{top}}^P(X, \{F_n\}).
\]    
Moreover, if $\{F_n\} $ is tempered, then     
\[ 
h_{\mathrm{top}}(X, G) = h_{\mathrm{top}}^B(X, \{F_n\}) = h_{\mathrm{top}}^P(X, \{F_n\}).
\]
\end{proposition}

The comparison among the three types of topological entropy for subsets is given as follows \cite{Dou2023ZhengZhou}.

\begin{proposition}\label{prop-entropy-B<P<UC}    
Let $(X,G)$ be a $G$-system and $\{F_n\} $ a F{\o}lner sequence of $G$. For any $K \subseteq X$:    
\begin{enumerate}        
    \item $h_{\mathrm{top}}^B(K, \{F_n\}) \le h_{\mathrm{top}}^P(K, \{F_n\})$.        
    
    \item If $\{F_n\} $ satisfies the growth condition $\lim\limits_{n \to \infty} \frac{|F_n|}{\log n} = \infty$, then 
    \[ 
    h_{\mathrm{top}}^P(K, \{F_n\}) \le h_{\mathrm{top}}(K, \{F_n\}).
    \]    
\end{enumerate}
\end{proposition}

\subsection{Hausdorff dimension}
Let $(Y,\rho)$ be a metric space and let $E\subseteq Y$.  For $s\ge0$ and
$\delta>0$, define
\[
\mathcal H_{\rho,\delta}^{s}(E)
=
\inf\left\{
\sum_{i=1}^{\infty}
\bigl(\operatorname{diam}_{\rho}U_i\bigr)^s:
E\subseteq\bigcup_{i=1}^{\infty}U_i,
\ \operatorname{diam}_{\rho}U_i\le\delta
\right\}.
\]
The \emph{$s$-dimensional Hausdorff measure} is
\[
\mathcal H_{\rho}^{s}(E)
=
\lim_{\delta\downarrow0}\mathcal H_{\rho,\delta}^{s}(E),
\]
and the \emph{Hausdorff dimension} of $E$ is
\[
\dim_H^{\rho}E
=
\inf\{s\ge0:\mathcal H_{\rho}^{s}(E)=0\}
=
\sup\{s\ge0:\mathcal H_{\rho}^{s}(E)=+\infty\}.
\]
We refer to \cite{Mattila1995,Falconer2014} for the standard theory.

\subsection{Lowering properties}\label{subsection-lowering-definitions}
Let $(X,G)$ be a $G$-system and let $\{F_n\} $ be a F{\o}lner
sequence.  We say that $(X,G)$ is
\begin{enumerate}
\item \emph{lowerable along $\{F_n\}$} if for every
$0\le h\le h_{\mathrm{top}}(X,G)$ there is a nonempty compact
$K_h\subseteq X$ with $h_{\mathrm{top}}(K_h,\{F_n\})=h$;
\item \emph{hereditarily lowerable along $\{F_n\}$} if for every nonempty
compact $K\subseteq X$ and every $0\le h\le h_{\mathrm{top}}(K,\{F_n\})$
there is a nonempty compact $K_h\subseteq K$ with
$h_{\mathrm{top}}(K_h,\{F_n\})=h$;
\item \emph{hereditarily uniformly lowerable along $\{F_n\}$} if the set
$K_h$ in the preceding item can always be chosen with at most one limit
point;
\item \emph{D-lowerable along $\{F_n\}$} if for every
$0\le h\le h_{\mathrm{top}}(X,G)$ there is a nonempty set $E_h\subseteq X$
with $h_{\mathrm{top}}^B(E_h,\{F_n\})=h$;
\item \emph{P-lowerable along $\{F_n\}$} if for every
$0\le h\le h_{\mathrm{top}}(X,G)$ there is a nonempty set $E_h\subseteq X$
with $h_{\mathrm{top}}^P(E_h,\{F_n\})=h$;
\item \emph{D-hereditarily lowerable along $\{F_n\}$} if for every
nonempty Souslin set $K\subseteq X$ and every
$0\le h\le h_{\mathrm{top}}^B(K,\{F_n\})$ there is a set $K_h\subseteq K$
with $h_{\mathrm{top}}^B(K_h,\{F_n\})=h$.
\end{enumerate}

Here and throughout, a Souslin set means an analytic subset of $X$. This class is well suited to the definition of D-hereditary lowerability. Souslin sets are closed under continuous images and preimages \cite[Section~14]{Kechris1995}, a useful property for transferring entropy-lowering results through factor maps. 
Howroyd's theorem \cite{Howroyd1995} guarantees that every Souslin set contains compact subsets of each prescribed Hausdorff dimension strictly below its own, a step that will be used later in the symbolic argument.
We do not introduce a uniform version of D-hereditary lowerability. Indeed, a compact set with at most one limit point is countable, and every countable subset has zero Bowen topological entropy \cite[Proposition~2.2]{Huang2014YeZhang}. Hence, the uniform requirement used for covering entropy cannot realize any positive value of Bowen topological entropy.

\subsection{Tail entropy and extension tools}\label{subsection-tail-extension-tools}
We collect here the notions and extension results needed for the last two sections. In particular, we recall conditional and tail topological entropy~\cite{Downarowicz2011}, relative topological entropy with its variational principles~\cite{Dooley2015Zhang,Yan2015}, and the quasi-symbolic extension characterization of asymptotic 
$h$-expansiveness~\cite{Downarowicz2023Zhang}. These tools will be applied first to 
D-hereditary lowerability for Bowen entropy, and then to hereditary uniform lowerability for compact subsets.

Let $(X,G)$ be a $G$-system. For subset $V \subseteq X$ and open covers $\mathcal{U,V} \in \mathcal C_X $, define $N( \mathcal U | \mathcal V) = \sup_{ V \in \mathcal V} N(\mathcal U,V) $. Let $\{F_n\}$ be a F{\o}lner sequence of $G$. One can verify that the function $F \mapsto \log N(\mathcal{U}_F|\mathcal{V}_F)$ is monotone and subadditive on $\mathcal{F}(G)$. Consequently, the limit
\[  h_{\mathrm{top}}(G,\mathcal{U}|\mathcal{V})= \lim_{n \to \infty} \frac{1}{|F_n|}\log N( \mathcal{U}_{F_n} | \mathcal{V}_{F_n} ) \]
exists. Note that $ h_{\mathrm{top}}(G,\mathcal{U}|\mathcal{V}) $ increases as $ \mathcal{U} $ is refined and decreases as $\mathcal{V}$ is refined.

The \emph{conditional topological entropy of $(X,G)$ with respect to $\mathcal{V}$} is defined as
\[ 
h_{\mathrm{top}}(G,X|\mathcal{V}) = \sup_{ \mathcal{U} } h_{\mathrm{top}}(G,\mathcal{U}|\mathcal{V}) ,
\]
where the supremum is taken over all finite open covers $\mathcal{U}$ of $X$. Observe that $ h_{\mathrm{top}} (X,G) = h_{\mathrm{top}}(G,X| \{X \}) $. 

The \emph{tail topological entropy of $(X,G)$} is defined by
\[ h_{\mathrm{top}}^* (X,G) = \inf_{\mathcal{V}} h_{\mathrm{top}}(G,X|\mathcal{V}) , \]
where the infimum is taken over all finite open covers $\mathcal{V}$ of $X$. The system $(X,G)$ is called \emph{asymptotically $h$-expansive} if $ h_{\mathrm{top}}^* (X,G) = 0 $.

Moreover, by the Ornstein-Weiss theorem \cite{OrnsteinWeiss1987,Gromov1999, Lindenstrauss2000Weiss}, the quantity $ h_{\mathrm{top}}(G,\mathcal{U}| \mathcal{V}) $ is independent of the choice of the F{\o}lner sequence $\{F_n\}$.

According to \cite{Dou2025WangZhang}, the tail topological entropy admits the following dynamical characterization.
\begin{lemma}
\label{prop-character-tailentropy}
Let $(X,G)$ be a $G$-system and $\{F_n\}$ be a F{\o}lner sequence for $G$. Then
\begin{align*}
h_{\mathrm{top}}^* (X,G) &= \lim_{ \varepsilon \to 0^+ } \sup_{x \in X} h_{\mathrm{top}} ( \Gamma_{\varepsilon, \{F_n\} } (x) , \{F_n\} ) \\
&= \lim_{ \varepsilon \to 0^+ } \sup_{x \in X} h_{\mathrm{top}} ( \Gamma_{\varepsilon } (x) , \{F_n\} ),
\end{align*}
where 
\[
\Gamma_{\varepsilon, \{F_n\}}(x) = \bigcap_{n \in \mathbb{N}} \bigcap_{g \in F_n} g^{-1} B(gx, \varepsilon)
\quad \text{and} \quad
\Gamma_{\varepsilon}(x) = \bigcap_{g \in G} g^{-1} B(gx, \varepsilon).   
\]
\end{lemma}

For a $G$-system $(Z,G)$, let $\mathcal{M}(Z)$ denote the space of Borel probability measures on $Z$, let $\mathcal{M}_G(Z)$ and $\mathcal{M}^e_G(Z)$ denote the sets of all $G$-invariant and ergodic $G$-invariant Borel probability measures on $Z$, respectively.
If $\pi\colon Y\to X$ is a factor map and $\nu\in\mathcal{M}(Y)$, we denote by $\pi_*\nu=\nu\circ\pi^{-1}$ the pushforward of $\nu$ under $\pi$.
An extension $\pi\colon (Y,G)\to (X,G)$ is said to be \emph{principal}, or equivalently $(X,G)$ is said to be a \emph{principal factor} of $(Y,G)$, if
$$
h_{\nu}(G,Y)=h_{\pi_*\nu}(G,X)
\quad\text{for every }\nu\in\mathcal{M}_G(Y).
$$
In this case, \(\pi\) is called a \emph{principal factor map}.

Now, fix a factor map $\pi: (Y,G) \to (X,G)$ and a F{\o}lner sequence $\{F_n\} $ for $G$. For an open cover $\mathcal{U} \in \mathcal{C}^o_Y$, define
\[
N(\mathcal{U}|\pi) = \max \left\{ N(\mathcal{U}, \pi^{-1}(x)) : x \in X \right\}.
\]
The \emph{relative topological entropy of $\mathcal{U}$ with respect to $\pi$} is defined as
\[
h_{\mathrm{top}}(G, \mathcal{U}|\pi) = \lim_{n \to \infty} \frac{1}{|F_n|} \log N(\mathcal{U}_{F_n}|\pi).
\]
By the Ornstein–Weiss theorem, this quantity is independent of the choice of the F{\o}lner sequence $\{F_n\} $.

The \emph{relative topological entropy of $(Y,G)$ with respect to $\pi$} and the \emph{conditional topological entropy of $(Y,G)$ with respect to $\pi$} are defined, respectively, as
\begin{align*}
h_{\mathrm{top}}(G,Y|\pi) &=\sup_{\mathcal{U}  } h_{\mathrm{top}}(G,\mathcal{U}|\pi) \\
h_{\mathrm{top}}^* (G,Y|\pi) &=\inf_{\mathcal{V} } \sup_{\mathcal{U}  } h_{\mathrm{top}}(G, \mathcal{U}|\pi^{-1}\mathcal{V} ) ,
\end{align*}
where $\mathcal{U}$ ranges over all finite open covers of $Y$ and $\mathcal{V}$ ranges over all finite open covers of $X$. 

The following result relates relative and conditional topological entropy \cite{Zhang2023Zhu}.

\begin{lemma}\label{lem-relative&conditional}
Let $\pi:(Y,G) \to (X,G)$ be a factor map between two $G$-systems. Then
\[ h_{\mathrm{top}}(G,Y|\pi) \le h_{\mathrm{top}}^*(G,Y|\pi) \le h_{\mathrm{top}}(G,Y|\pi) + h_{\mathrm{top}}^*(X,G) . \]
In particular, if $(X,G)$ is asymptotically $h$-expansive, then 
\[ h_{\mathrm{top}}(G,Y|\pi) = h_{\mathrm{top}}^*(G,Y|\pi) .\]
\end{lemma}

We will also use the following results from \cite{Zhu2021}.

\begin{lemma}\label{lem-condition-asyhex}
Let $\pi:(Y,G) \to (X,G)$ be a factor map between two $G$-systems. Then
\begin{enumerate}
\item $h_{\mathrm{top}}^* (Y,G) \le h_{\mathrm{top}}^* (G,Y|\pi)$.
\item $h_{\mathrm{top}}^* (X,G) \le h_{\mathrm{top}}^*(Y,G) + h_{\mathrm{top}}^*(G,Y |\pi) $.
\end{enumerate}
\end{lemma}

Given $\nu \in \mathcal{M}_G(Y)$, the preimage $\pi^{-1} \mathcal{B}_X$ is clearly a $G$-invariant sub-$\sigma$-algebra of $\mathcal{B}_Y$. For a finite measurable partition $\beta$ of $Y$, the \emph{measure-theoretic conditional $\nu$-entropy of $\beta$ with respect to $\pi$} is defined by
\[ 
h_{\nu}(G,\beta|\pi):= h_{\nu}(G,\beta| \pi^{-1} \mathcal{B}_X ) =\lim\limits_{n \to \infty} \frac{1}{|F_n|} H_{\nu}( \beta_{F_n}| \pi^{-1} \mathcal{B}_X ). 
\]
The \emph{measure-theoretic conditional $\nu$-entropy of $(Y,G)$ with respect to $\pi$} is then defined as 
\[ 
h_{\nu}(G,Y|\pi)=\sup\limits_{\beta} h_{\nu}(G,\beta|\pi), 
\]
where the supremum is taken over all finite measurable partitions $\beta$ of $Y$.

For a fixed F{\o}lner sequence $\{F_n\}$ of $G$, the function $ x \in X \mapsto h_{\mathrm{top}}( \pi^{-1}(x) , \{F_n\} )  $ is measurable on $X$. 
The variational principle for the conditional entropy mentioned below was essentially established by Dooley and Zhang in their work \cite{Dooley2015Zhang}. Yan \cite{Yan2015} also subsequently provided this result.

\begin{proposition}\label{prop-VP-yan}
Let $\pi: (Y,G) \to (X,G) $ be a factor map between two $G$-systems and $\{F_n\} $ a F{\o}lner sequence for $G$ with $e_G \in F_1 \subseteq F_2 \subseteq \cdots$. Then
\begin{enumerate}
\item $h_{\mathrm{top}}(G,Y|\pi) =\sup\limits_{\nu \in \mathcal{M}_G(Y)} h_{\nu}(G,Y|\pi)$.
\item For any $\mu \in \mathcal{M}_G(X)$,  
\[ 
\sup \{ h_{\nu}(G,Y|\pi) : \nu \in \mathcal{M}_G(Y),\, \pi \nu = \mu \} =\int_{X} h_{\mathrm{top}}(\pi^{-1}(x),\{F_n\}) d\mu (x) .
\]
Equivalently, by the classical variational principle,
\[  
\sup \{ h_{\nu}(G, Y) : \nu \in \mathcal{M}_G(Y),\, \pi \nu = \mu \} = h_{\mu}(G,X) + \int_{X} h_{\mathrm{top}}(\pi^{-1}(x),\{F_n\}) d\mu (x). 
\]
\item If $\{F_n\} $ is also tempered, then
\[ 
h_{\mathrm{top}}(G,Y|\pi) =\sup\limits_{x \in X} h_{\mathrm{top}} (\pi^{-1}(x) , \{F_n\} ) = \sup\limits_{ \mu \in \mathcal{M}_G (X) } \int _X h_{\mathrm{top}} (\pi^{-1}(x),\{F_n\} ) d \mu(x) .
\]
\end{enumerate}
\end{proposition}

\begin{definition}
A topological joining $\overline{Y} = Y \vee \mathbf{T}$ of a subshift $Y$ with a tiling system $\mathbf{T}$ of zero entropy is called a \emph{quasi-symbolic system}. A topological extension that is itself a quasi-symbolic system is called a \emph{quasi-symbolic extension}.
\end{definition}

The following equivalent characterization of asymptotic $h$-expansiveness was established by Downarowicz and Zhang \cite{Downarowicz2023Zhang}.

\begin{lemma}\label{lem-equicondition-asyhexpan}
A $G$-system $(X,G)$ is asymptotically $h$-expansive if and only if it admits a principal quasi-symbolic extension.
\end{lemma}

\section{Lowering entropy for finite-entropy and expansive $G$-systems}\label{section-Finite-entropy-system}
In this section, we establish the direct entropy-lowering results of the paper. We first prove the entropy distribution principles and relative entropy estimates needed for systems of finite topological entropy, and then prove Theorem~\ref{thm-main-finite}. We next consider expansive actions and prove hereditary uniform lowerability by constructing, inside a given compact set, subsets of prescribed entropy that are countable and have at most one accumulation point. These results also provide the basic constructions for the study of hereditary lowering and tail entropy in the subsequent sections.

\subsection{Distribution principles}\label{subsection-bridge}
The following result gives the quantitative relation between the topological entropy and Bowen topological entropy of subsets.
\begin{lemma}\label{lem3.1}
Let $(X,G) $ be a $G$-system, $\{ F_n\} $ be a F{\o}lner sequence of $G$, $K \subseteq X$ and $\mathcal{U} \in \mathcal{C}_X^o$. Then 
\[
h_{\mathrm{top}}^B ( K , \mathcal{U} , \{F_n\} ) \leq \liminf\limits_{n \to \infty} \frac{1}{|F_n|} \log N ( \mathcal{U}_{F_n} , K ) \leq h_{\mathrm{top}} ( K , \mathcal{U} , \{F_n\} ).     
\]
\end{lemma}
\begin{proof}
If $K = \emptyset$, this is obvious. Now we assume $ K \neq \emptyset $. For every $ n \in \mathbb{N} $, set 
\[
\mathcal{T}_n = \{ A_1, A_2, \cdots, A_{N(\mathcal{U}_{F_n} , K ) } \} \subseteq \mathcal{U}_{F_n}     
\] 
with $K \subseteq \bigcup \mathcal{T}_n$. For each $A \in \mathcal{T}_n$, there is a $\mathbb{U}_A \in \mathcal{W}_{F_n} (\mathcal{U})$ such that 
\[
A = X( \mathbb{U}_A ) \text{ and } m( \mathbb{U} _A ) = |F_n| . 
\]
Thus $\Lambda = \{ \mathbb{U}_A : A \in \mathcal{T}_n \} \subseteq \bigcup_{j \geq n} \mathcal{W}_{F_j} (\mathcal{U})$ covers $K$. For each $s \geq 0$, one has 
\begin{align}
\mathcal{M} (K, \mathcal{U}, n, s, \{F_n\}) &\leq \sum_{\mathbb{U}_A \in \Lambda} e^{-s m(\mathbb{U}_A )} \notag \\ &= \sum_{ A \in \mathcal{T}_n } e^{-s |F_n|} = N(\mathcal{U}_{F_n} , K) e^{-s |F_n|}, \notag 
\end{align}
so 
\[
\mathcal{M} ( K, \mathcal{U}, s, \{F_n\} ) \leq \liminf\limits_{n \to \infty } e^{-|F_n|[s- \frac{1}{|F_n|} \log N(\mathcal{U}_{F_n} , K) ]} .   
\]
If $s > \liminf\limits_{n \to \infty} \frac{1}{|F_n|} \log N(\mathcal{U}_{F_n},K)$, then $\mathcal{M}(K, \mathcal{U}, s, \{F_n\}) = 0$. Therefore 
\[
h_{\mathrm{top}}^B(K, \mathcal{U}, \{F_n\}) \leq \liminf\limits_{n \to \infty} \frac{1}{|F_n|} \log N(\mathcal{U}_{F_n} , K) ,    
\] 
which ends the proof.
\end{proof}

Let $(X,G)$ be a $G$-system and $\{F_n\}_{n \in \mathbb N } $ be a F{\o}lner sequence of $G$ with $e_G \in F_1 \subseteq F_2 \subseteq \cdots$. The following lemmas are useful for proving Theorem~\ref{thm-main-finite} in this section.

\begin{lemma}\label{lem-nonuniform-massdis-principle}
Let $d>0$, $M \in \mathbb{N}$, $K \subseteq X$, $\alpha \in \mathcal{P}_X$, $\mathcal{U} \in \mathcal{C}_X$ and $\mu \in \mathcal{M}(X)$. Suppose that each element of $\mathcal{U}$ has a non-empty intersection with at most $M$ elements of $\alpha$. If there exists $K_{\mu} \subseteq K$ such that $K_{\mu}$ has a positive outer $\mu$-measure (that is, $\mu^* (K_{\mu}) > 0$) and for each $x \in K_{\mu}$ there is a constant $c(x) > 0$ such that 
\[
\mu ( \alpha_{F_n} (x) ) \leq c(x)\cdot e^{-|F_n|d}, \text{ } \forall n \in \mathbb{N} .    
\]
Then $h_{\mathrm{top}}^B (K, \mathcal{U}, \{F_n\}) \geq d - \log M$. In particular, $h_{\mathrm{top}}^B (K,\alpha, \{F_n\}) \geq d$. 
\end{lemma}

\begin{proof}
Assuming $d - \log M > 0$ and for each $k \in \mathbb{N}$, We set 
\[
K_{\mu}^k = \{ x \in K_{\mu} : c(x) \leq |F_k| \}.    
\]
Thus, $K_{\mu}^1 \subseteq K_{\mu}^2 \subseteq \cdots$ and $K_{\mu} = \bigcup_{k\in \mathbb{N}} K_{\mu}^k$. By $ \mu ^* (K_{\mu}) > 0$, there exists $N \in \mathbb{N}$ with $\mu ^*(K_{\mu}^N) > 0 $. 
		
For any $k \in \mathbb{N}$ and $\Lambda \subseteq \bigcup_{j \geq k} \mathcal{W}_{F_j}(\mathcal{U})$ covers $K$. 
We denote the index set of $\mathbb{U} = \{U_g\} _{g \in F_n}$ to be $F_{n(\mathbb{U})}$ in convention. 
If $\mathbb{U} \in \Lambda$ satisfied $X(\mathbb{U}) \cap K_{\mu}^N \neq \emptyset$, then for each $l \in \mathbb{N}$ and $B \in \alpha_{F_{\min\{l, n(\mathbb{U})  \}}}$ with $B \cap (X(\mathbb{U}) \cap K_{\mu}^N) \neq \emptyset$, one can choose 
\[
x_B \in B \cap ( X(\mathbb{U}) \cap K_{\mu}^N )     
\]
and
\begin{align}
\mu(B) &= \mu (\alpha_{F_{\min\{l, n(\mathbb{U})\}}} (x_{B})) \notag \\ 
&\leq c(x_B)\cdot e^{- \min \{ |F_l| , |F_{n( \mathbb{U} )}| \} \cdot d} \notag \\ 
&\leq |F_N| e^{- \min \{ |F_l| , |F_{n( \mathbb{U})}| \}  \cdot d} \notag
\end{align}
For $\mathbb{U} \in \Lambda \subseteq \bigcup_{j \geq k} \mathcal{W}_{F_j} (\mathcal{U})$, that is $X(\mathbb{U}) \in \mathcal{U}_{F_{n(\mathbb{U})}}$ for some $n(\mathbb{U}) \geq k$ and $|F_{n(\mathbb{U})}| = m(\mathbb{U})$, one has
\[
\left | \{ V \in \alpha_{ F_{ min\{l, n(\mathbb{U})  \} } } : V \cap ( X(\mathbb{U}) \cap K_{\mu}^N ) \neq \emptyset \} \right | \leq M^{ \min \{ |F_l|, |F_{ n( \mathbb{U} ) }| \} } = M^{ \min \{ |F_l|, m( \mathbb{U} ) \} }  ,   
\]
so  
\begin{align}
\mu^* (X(\mathbb{U}) \cap K_{\mu}^N) &\leq M^{\min \{ |F_l|, m(\mathbb{U})\}} \cdot |F_N| e^{- \min \{|F_l|, m( \mathbb{U})\} \cdot d} \notag \\ 
&= |F_N| e^{- \min \{|F_l|, m(\mathbb{U})\} \cdot (d - \log M)} .\notag
\end{align}
For the arbitrariness of $l \in \mathbb{N}$, we obtain that $\mu^* (X(\mathbb{U}) \cap K_{\mu}^N) \leq |F_N| e^{ -(d - \log M) m(\mathbb{U})}$ for each $\mathbb{U} \in \Lambda$ which satisfies $X(\mathbb{U}) \cap K_{\mu}^N$. Furthermore, 
\begin{align}
\sum_{\mathbb{U} \in \Lambda} e^{-(d - \log M) m(\mathbb{U})} &\geq \sum_{\mathbb{U} \in \Lambda : X(\mathbb{U}) \cap K_{\mu}^N \neq \emptyset} e^{-(d - \log M) m(\mathbb{U})} \notag \\ 
&\geq \sum_{\mathbb{U} \in \Lambda : X(\mathbb{U}) \cap K_{\mu}^N \neq \emptyset} \frac{1}{|F_N|} \mu^* ( X(\mathbb{U}) \cap K_{\mu}^N) \notag \\ 
&\geq \frac{1}{|F_N|} \mu^* (K_{\mu}^N). \notag 
\end{align}
By the choice of $\Lambda$, for every $k \in \mathbb{N} $, one has 
\[
\mathcal{M} (K, \mathcal{U}, k, d- \log M, \{F_n\}) \geq \frac{1}{|F_N|} \mu^* (K_{\mu}^N) > 0 ,    
\]
thus 
\[
\mathcal{M} (K, \mathcal{U}, d- \log M, \{F_n\} ) \geq \frac{1}{|F_N|} \mu^* (K_{\mu}^N) > 0 .   
\]
Therefore, 
\[
h_{\mathrm{top}}^B (K, \mathcal{U}, \{F_n\}) \geq d - \log M .   
\] 
In particular, setting $\mathcal{U}=\alpha$ yields $h_{\mathrm{top}}^B(K,\alpha,\{F_n\}) \ge d$. This completes the proof.
\end{proof}

\begin{lemma}\label{lem-uniform-massdis-principle}
Let $c>0$, $d>0$, $K \subseteq X$ and $\alpha \in \mathcal{P}_X$. If there exists $\mu \in \mathcal{M}(X)$ such that 
\[
\mu ( \alpha_{F_n} (x) ) \geq c \cdot e^{ -|F_n|d } \text{, for any } x \in K \text{ and } n \in \mathbb{N},  
\]
then $h_{\mathrm{top}} (K, \alpha, \{F_n\}) \leq d$.
\end{lemma}

\begin{proof}
Assume that $K \neq \emptyset$, for every $n \in \mathbb{N}$, we set $k = N (\alpha_{F_n}, K)$ and $\mathcal{T}_n = \{ A_1, \cdots, A_k\} \subseteq \alpha_{F_n}$ with $K \subseteq \bigcup \mathcal{T}_n$. Then $A_i \cap K \neq \emptyset$, $1 \leq i \leq k$.
		
For $i \in \{1,\cdots,k\}$, one can choose $x_i \in A_i \cap K$, then 
\[
\mu(A_i) = \mu(\alpha_{F_n} (x_i)) \geq c \cdot e^{-|F_n|d}.   
\]
So $1 \geq \sum_{i=1}^{k} \mu(A_i) \geq N(\alpha_{F_n}, K) \cdot c e^{-|F_n|d}$. And hence $N(\alpha_{F_n}, K) \leq \frac{e^{|F_n|d}}{c}$.
Therefore, we obtain 
\[
h_{\mathrm{top}}(K, \alpha, \{F_n\}) = \limsup\limits_{n \to \infty} \frac{1}{|F_n|} \log N( \alpha_{F_n}, K ) \leq d.   
\]
\end{proof}

\begin{remark}\label{rmk3.1}
The next subsection relies on the mass distribution principles—namely, the non-uniform and uniform versions formulated as Lemma \ref{lem-nonuniform-massdis-principle} and Lemma \ref{lem-uniform-massdis-principle}, respectively. These principles link the lowerable property to ergodic theory and are essentially established in \cite{Pesin1997}.
\end{remark}

\subsection{Lowering entropy for finite-entropy systems}\label{subsection-low-Dlow-Plow}
Consider a \textit{measure-theoretic dynamical system} $(X, \mathcal{B}, \mu, G)$, where $(X, \mathcal{B}, \mu)$ is a probability space and $G$ is a countably infinite discrete amenable group. The system is termed \textit{measure-preserving} under the action of $G$, meaning that every group element acts as a measure-preserving transformation on the space.

Let $(X, G)$ be a topological $G$-system.  
Given any $\mu \in \mathcal{M}_G(X)$, the quadruple $(X, \mathcal{B}_X, \mu, G)$ forms a measure-preserving dynamical system.
By the amenability of $G$, the set $\mathcal{M}_G^e(X)$ is nonempty. Furthermore, both $\mathcal{M}_G(X)$ and $\mathcal{M}(X)$ are convex, compact, and metrizable spaces when equipped with the weak $*$ topology.
Let $\{F_n\}$ be a F{\o}lner sequence of $G$.
	
Given $ \alpha \in \mathcal{P}_X $, $ \mu \in \mathcal{M}_G(X) $ and a sub-$\sigma$-algebra $\mathcal{A} \subseteq \mathcal{B}_{\mu}$, where $ \mathcal{B}_{\mu} $ is the completion of $ \mathcal{B}_X $ under $ \mu $, the {\it conditional information function of $ \alpha $ relevant to $ \mathcal{A} $} is defined by 
$$ I_{\mu} (\alpha | \mathcal{A}) (x) = \sum_{A \in \alpha} - 1_A (x) \log \mathbb{E}_{\mu} ( 1_A| \mathcal{A} ) (x) ,  $$
where $ \mathbb{E}_{\mu} ( 1_A| \mathcal{A} ) $ is the expectation of $ 1_A $ with respect to $ \mathcal{A} $. 
Let $$ H_{\mu} ( \alpha| \mathcal{A} ) = \int_{X} I_{\mu} (\alpha | \mathcal{A}) (x) d \mu(x) = \sum_{A \in \alpha} \int_{X} - \mathbb{E}_{\mu} ( 1_A| \mathcal{A} ) \log \mathbb{E}_{\mu} ( 1_A| \mathcal{A} ) d \mu.$$ 
Observe that $H_{\mu}(\alpha | \mathcal{A})$ is monotonic in its arguments: it is non-decreasing in $\alpha$ (with respect to refinement of the partition) and non-decreasing in $\mathcal{A}$ (with respect to coarsening of the $\sigma$-algebra). 
The {\it measure-theoretic $ \mu $-entropy of $\alpha$ relevant to $\mathcal{A}$} is defined as 
\[ h_{\mu} (G,\alpha| \mathcal{A} ) = \lim\limits_{n \to \infty } \frac{1}{|F_n|} H_{\mu} ( \alpha_{F_n} | \mathcal{A} ) = \inf\limits_{F \in \mathcal{F}(G)} \frac{1}{|F|} H_{\mu} (\alpha_{F} | \mathcal{A}) .  \]
Since the map
$F \mapsto H_{\mu}(\alpha_F \mid \mathcal{A})$
is monotone, nonnegative, $G$-invariant, and subadditive on $\mathcal{F}(G)$, the Ornstein--Weiss lemma implies that the limit exists and equals the corresponding infimum; see \cite{OrnsteinWeiss1987,Huang2011YeZhang}. 
The {\it measure-theoretic $ \mu $-entropy of $ (X,G) $ relevant to $ \mathcal{A} $} is defined as
\[ h_{\mu}(G,X|\mathcal{A}) = \sup\limits_{\alpha \in \mathcal{P}_X } h_{\mu} (G,\alpha| \mathcal{A} )  .\]
	
Set 
\[ H_{\mu} ( \mathcal{U}| \mathcal{A} ) = \inf\limits_{\alpha \in \mathcal{P}_X : \alpha \succeq \mathcal{U}} H_{\mu} ( \alpha| \mathcal{A} ) \]
for $ \mathcal{U} \in \mathcal{C}_X $. Clearly, $ H_{\mu} (\mathcal{U} | \mathcal{A}) $ is non-decreasing in $ \mathcal{U} $ and non-increasing in $ \mathcal{A} $. 
Similarly, the {\it measure-theoretic $\mu^{-}$-entropy of $\mathcal{U}$ relevant to $\mathcal{A}$} is defined by \[ h_{\mu}^{-} (G,\mathcal{U} |\mathcal{A} ) = \lim\limits_{n \to \infty} \frac{1}{|F_n|} H_{\mu} (\mathcal{U}_{F_n} | \mathcal{A} ) = \inf\limits_{F \in \mathcal{F}(G)} \frac{1}{|F|} H_{\mu} (\mathcal{U}_{F} | \mathcal{A} ) \]
and this value is independent of the choice of F{\o}lner sequence $ \{F_n\}_{n \in \mathbb N} $. As shown in \cite[Theorem 3.5]{Huang2011YeZhang}, it follows that \[ h_{\mu} (G,X | \mathcal{A}) = \sup\limits_{\mathcal{U} \in \mathcal{C}_X^o} h_{\mu}^{-} (G, \mathcal{U} | \mathcal{A}) .\]
Consequently, if $ \{ \mathcal{U}_m\}_{m \in \mathbb N} \subseteq \mathcal{C}_X^o $ satisfies $ \lim\limits_{m \to \infty} \| \mathcal{U}_m \| = 0 $, then 
\begin{align}\label{measure-entropy-cover}
h_{\mu} (G,X | \mathcal{A}) = \lim\limits_{n \to \infty} h_{\mu}^{-} (G, \mathcal{U}_m | \mathcal{A}) .
\end{align}

When $\mathcal{A} = \{X, \emptyset\} \pmod{\mu}$, we omit $\mathcal{A}$ in the notation and write $H_{\mu}(\mathcal{U})$, $h_{\mu}(G, \alpha)$, $h_{\mu}^{-}(G, \mathcal{U})$, and $h_{\mu}(G, X)$ for $H_{\mu}(\mathcal{U} | \mathcal{A})$, $h_{\mu}(G, \alpha | \mathcal{A})$, $h_{\mu}^{-}(G, \mathcal{U} | \mathcal{A})$, and $h_{\mu}(G, X | \mathcal{A})$, respectively.

Following the proof technique of \cite[Theorem 15.11]{Glasner2003} in a similar manner, we proceed to demonstrate the following pivotal result.
\begin{proposition}\label{prop-exist-factor}
Let $(X, \mathcal{B}, \mu, G)$ be an ergodic measure-preserving dynamical system with $G$ amenable. Then for every $s \in [0, h_{\mu}(G, X)]$, there exists a factor $(Z,\mathcal{Z},\nu,G)$ of $(X, \mathcal{B}, \mu, G)$ such that $h_{\nu}(G, Z) = s$.
\end{proposition}
\begin{proof}
Fix an integer $l \ge 2$ and denote by $\mathcal{P}_l$ the set of all measurable partitions of $X$ consisting of $l$ sets. For $ \alpha = \{A_1,\dots,A_l\} $, $\beta=\{B_1,\dots,B_l\} \in \mathcal{P}_l$, we define 
\[ 
d_{\text{part} }(\alpha,\beta) = \mu ( \alpha \bigtriangleup \beta ) = \frac{1}{2} \sum_{j=1}^l \mu(A_j \bigtriangleup B_j) .
\]
Then $(\mathcal{P}_l , d_{\text{part} })$ forms a complete metric space. We also set 
\[ 
d_{\text{ent}}(\alpha,\beta) = H_{\mu}(\alpha | \beta) +H_{\mu}(\beta | \alpha), 
\]
which defines a metric on tha space of all (finite or countable) partitions with finite entropy. 
Given any $\varepsilon > 0$, there exists a $\delta=\delta(\varepsilon,l)>0$ such that for any pair $\alpha,\beta \in \mathcal{P}_l$ with $d_{\text{part} }(\alpha,\beta)<\delta$, one has $ d_{\text{ent} }(\alpha,\beta)< \varepsilon $.
Consequently, the function $\alpha \mapsto h_{\mu}(G,\alpha) $ is continuous on $\mathcal{P}_l$ with respect to the metric $d_{\text{part}}$. 

By Kreiger's finite generator theorem for amenable group actions \cite[Theorem 1.1, 2.3]{Seward2019}, we may assume that $ h_{\mu}(G,X)=h_{\mu}(G,\alpha) $ for some finite partition $\alpha$. 
Let $l = |\alpha|$.
In the metric space $( \mathcal{P}_l , d_{\text{part} } )$, there exists a path $ \alpha_{t} $, $0 \le t \le 1$ with $\alpha_0 = \{X, \emptyset, \dots, \emptyset\}$ and $\alpha_1 = \alpha$. 

Since the function $\alpha \mapsto h_{\mu}(G, \alpha)$ is continuous with respect to $d_{\text{part}}$, the function $f(t) := h_{\mu}(G, \alpha_t)$ is continuous on $[0, 1]$. Note that $f(0) = 0$ and $f(1) = h_{\mu}(G, X)$. By the Intermediate Value Theorem, the image of $f$ contains the entire interval $[0, h_{\mu}(G, X)]$. Therefore, for any $s \in [0, h_{\mu}(G, X)]$, there exists $t_0 \in [0, 1]$ such that $f(t_0) = s$.
Let $(Z,\mathcal{Z},\nu,G)$ be the factor induced by the $G$-invariant $\sigma$-algebra $(\alpha_{t_0})_G = \bigvee_{g \in G} g\alpha_{t_0}$. Then $h_{\nu}(G, Z) = s$, which completes the proof.
\end{proof}
	
Let $\mu \in \mathcal{M}_G(X)$. Then $ (X, \mathcal{B}_{\mu}, \mu, G) $ forms a Lebesgue system. Consider a countable family $\{\alpha_i\}_{i \in I}$ in $\mathcal{P}_X$. The partition $ \alpha = \bigvee_{i \in I} \alpha_i $, defined by
\[ 
\alpha = \bigl\{ \cap_{i\in I} A_i : A_i \in \alpha_i \text{ for each } i \in I \bigr\},
\] 
is called a {\it measurable partition}. 
The subsets $ A \in \mathcal{B}_{\mu} $ that are unions of atoms of $ \alpha $ constitute a sub-$ \sigma $-algebra of $ \mathcal{B}_{\mu} $, denoted by $ \widehat{\alpha} $ (or simply $ \alpha $ if no ambiguity arises). 
It is a standard result that every sub-$\sigma$-algebra of $\mathcal{B}_{\mu}$ coincides, modulo sets of measure zero, with one constructed in this manner.
	
Let $\mathcal{A}$ be a sub-$\sigma$-algebra of $\mathcal{B}_{\mu}$, and let $\alpha$ be a measurable partition of $X$ such that $ \widehat{\alpha} = \mathcal{A} \pmod{\mu}$. 
Consider the disintegration $ \mu = \int_{X} \mu_{x} d \mu (x) $ with respect to $ \mathcal{A} $, where $\mu_{x} \in \mathcal{M}(X) $ and $ \mu_{x}(\alpha(x)) $ for $ \mu $-a.e. $ x \in X $. This disintegration is characterized by the following properties:
\begin{enumerate}
    \item[1.] For every $f \in L^1(X, \mathcal{B}_X, \mu)$, we have $f \in L^1(X, \mathcal{B}_X, \mu_x)$ for $\mu$-a.e. $x \in X$, and the map $x \mapsto \int_X f(y) d\mu_x(y)$ belongs to $L^1(X, \mathcal{A}, \mu)$.
    \item[2.] For every $f \in L^1(X, \mathcal{B}_X, \mu)$, the conditional expectation satisfies
    \[
        \mathbb{E}_{\mu} (f | \mathcal{A})(x)=\int_{X} f d \mu_{x} \quad \text{for } \mu \text{-a.e. } x \in X.
    \]
\end{enumerate}
Consequently, for any $ f \in L^1(X,\mathcal{B}_X,\mu) $, \[
\int_{X} \left( \int_{X} f d \mu_{x} \right) d \mu (x) = \int_{X} f d \mu .     
\] 
For $ \mu $-a.e. $ x \in X $, define the set $ \Gamma_{x} = \{ y \in X : \mu_y=\mu_x \} $. Then $ \mu_{x}(\Gamma_{x}) = 1$ for $\mu$-a.e. $ x \in X $. Hence for every $ f \in L^1(X,\mathcal{B}_X,\mu) $ and $ \mu $-a.e. $ x \in X $, we have
\begin{align}\label{Exp-constant-fiber}
	\mathbb{E}_{\mu} ( f | \mathcal{A} ) (y) = \int_{X} f d \mu_y = \int_{X} f d \mu_x = \mathbb{E}_{\mu} ( f | \mathcal{A} )(x) \quad \text{for } \mu_x\text{-a.e. }y \in X.
\end{align}
In particular, if $ f $ is $ \mathcal{A} $-measurable, then for $\mu$-a.e. $ x \in X $, $ f(y) = f(x) $ for $ \mu_x $-a.e. $ y \in X $.
	
Denote by $I_{\mu}$ the family of sets
\[
I_{\mu} = \{ A \in \mathcal{B}_{\mu} : \mu( gA \bigtriangleup A ) =0 \text{ for all } g \in G \} .    
\] 
It is clear that $I_{\mu}$ is a $G$-invariant sub-$\sigma$-algebra of $\mathcal{B}_{\mu}$. Let $ \alpha $ be a measurable partition of $ X $ such that $ \widehat{\alpha} = I_{\mu}  \pmod{\mu} $. 
It is easy to see that $\alpha $ is $G$-invariant, i.e., $ g \alpha = \alpha $ for all $ g \in G $. Consider the disintegration of $\mu$ with respect to $I_{\mu}$, $ \mu = \int_{X} \mu_{x} d \mu (x) $, where $ \mu_{x} \in \mathcal{M}_G^e(X) $ for $\mu$-a.e. $x \in X$. This disintegration is called the {\it ergodic decomposition} of $\mu$ (see, for example, \cite[Theorem 3.22]{Glasner2003}).
	
The following presents the relative versions of the Shannon–McMillan–Breiman theorem and the Brin–Katok entropy formula for countable discrete amenable group actions \cite{Lindenstrauss2001, Weiss2003,Dou2026WangZhou}.
\begin{theorem}\label{thm-re-SMB}
Let $(X,G) $ be a $G$-system, $\mu \in \mathcal{M}_G(X)$, $\{F_n\} $ a tempered F{\o}lner sequence of $G$ satisfying the growth condition 
\[
\lim\limits_{n \to \infty } \frac{|F_n|}{\log n} = \infty. 
\] 
Let $ \alpha \in \mathcal{P}_X $ and $ \mathcal{A} \subseteq \mathcal{B}_{\mu} $ be a $G$-invariant sub-$\sigma$-algebra (i.e. $g^{-1} \mathcal{A} =\mathcal{A} \text{ for all }g \in G $). Then there exists a $G$-invariant function $ f \in L^1(\mu) $ with $ \int_{X} f d \mu = h_{\mu} (G, \alpha | \mathcal{A}) $, such that 
\[
\lim\limits_{n \to \infty} \frac{1}{|F_n|} I_{\mu} (\alpha_{F_n} | \mathcal{A}) (x) = f(x)     
\]
for $\mu$-a.e. $x \in X$ and in $L^1(\mu)$. Moreover, if $\mu$ is ergodic then 
\[
\lim\limits_{n \to \infty}\frac{1}{|F_n|} I_{\mu}(\alpha_{F_n} | \mathcal{A}) (x) = h_{\mu}(G,\alpha|\mathcal{A}) \quad \text{for }\mu \text{-a.e. }x \in X.    
\] 
\end{theorem}

\begin{theorem}\label{thm-re-BrinKatok}
Let $(X,G) $ be a $G$-system, $\mu \in \mathcal{M}^e_G (X)$, $\{F_n\} $ a tempered F{\o}lner sequence of $G$ satisfying the growth condition 
\[
\lim\limits_{n \to \infty } \frac{|F_n|}{\log n} = \infty.    
\]
Let $ \mathcal{A} \subseteq \mathcal{B}_{\mu} $ be a $G$-invariant sub-$\sigma$-algebra, and let 
\[ 
\mu=\int_X\mu_x d\mu(x)
\]
be the disintegration of $\mu$ over $\mathcal{A}$. 
Then for $\mu$-a.e. $x\in X$, 
\[ 
\lim_{\varepsilon\to 0}\limsup_{n\to \infty} -\frac{1}{|F_n|}\log \mu_x (B_{F_n}(x,\varepsilon) ) = \lim_{\varepsilon\to 0}\liminf_{n\to \infty} -\frac{1}{|F_n|}\log \mu_x (B_{F_n}(x,\varepsilon) ) = h_{\mu}(G,X|\mathcal{A} ) .
\]
\end{theorem}

By combining the result of Theorem \ref{thm-re-SMB} with the proof strategy of Theorem \ref{thm-re-BrinKatok}, one obtains a non-ergodic version of the relative Brin–Katok entropy formula presented below.
\begin{theorem}
Let $(X,G) $ be a $G$-system, $\mu \in \mathcal{M}_G (X)$, $\{F_n\} $ a tempered F{\o}lner sequence of $G$ satisfying the growth condition 
\(
\lim\limits_{n \to \infty } \frac{|F_n|}{\log n} = \infty.    
\)
Let $ \mathcal{A} \subseteq \mathcal{B}_{\mu} $ be a $G$-invariant sub-$\sigma$-algebra, and let 
\( 
\mu=\int_X\mu_x d\mu(x)
\)
be the disintegration of $\mu$ over $\mathcal{A}$. 
Then 
\[
\underline{h}_{\mu_x } (X,\{F_n\}) = \overline{h}_{\mu_x  } (X,\{F_n\})
\]
for $\mu$-a.e. $x \in X$ and $\int_X \underline{h}_{\mu_x } (X,\{F_n\}) d \mu(x) = \int_X \overline{h}_{\mu_x} (X,\{F_n\}) d \mu(x) =  h_{\mu}(G,X|\mathcal{A} ) $, where
\begin{align*}
\underline{h}_{\mu_x} (X,\{F_n\}) &= \int_X \lim_{\varepsilon\to 0}\limsup_{n\to \infty} -\frac{1}{|F_n|}\log \mu_x (B_{F_n}(y,\varepsilon) ) d \mu_x (y), \\ \overline{h}_{\mu_x } (X,\{F_n\}) &= \int_X \lim_{\varepsilon\to 0}\liminf_{n\to \infty} -\frac{1}{|F_n|}\log \mu_x (B_{F_n}(y,\varepsilon) ) d \mu_x (y).    
\end{align*} 
\end{theorem}

\begin{remark}
Under the assumptions of Theorem \ref{thm-re-SMB}, let $ \mu = \int_{X} \mu_{x}^e d \mu (x) $ be the ergodic decomposition of $\mu $. The $G$-invariant function $f$ obtained in the theorem admits the explicit expression
\[
f(x) = h_{\mu_{x}^e} (G, \alpha | \mathcal{A}) .    
\]
In other words, for $ \mu$-almost every $x$, the limit of the conditional information function equals the conditional measure-theoretic entropy with respect to the ergodic component $\mu_x^e$.
\end{remark}
	
\begin{proposition}\label{prop-lowerable-essential}
Let $(X,G) $ be a $G$-system and $ \mu \in \mathcal{M}_G(X) $. Let $\{F_n\} $ be a tempered F{\o}lner sequence of $G$ such that $e_G \in F_1 \subseteq F_2 \subseteq \cdots $ and $ \lim\limits_{n \to \infty } \frac{|F_n|}{\log n} = \infty $. Suppose $\mathcal{A} \subseteq \mathcal{B}_{\mu} $ is a $G$-invariant sub-$\sigma$-algebra and let $ \mu = \int_{X} \mu_{x} d \mu (x) $ be the disintegration of $ \mu $ over $ \mathcal{A} $. Then the following hold:
\begin{enumerate}
\item[(1)] Let $ \mathcal{U} \in \mathcal{C}_X $ and $ \alpha \in \mathcal{P}_X $ such that each element of $ \mathcal{U} $ intersects at most $ M $ elements of $ \alpha $ for some $M \in \mathbb{N}$. 
Let $ f^{\{F_n\}}_{ \alpha , \mathcal{A} }(x) $ be the function from Theorem \ref{thm-re-SMB} for $ \{F_n\} $, $ \alpha $ and $ \mathcal{A} $. Then for $ \mu $-a.e. $ x \in X $ and any $ K_x \in \mathcal{B}_X $ with $ \mu_{x}(K_x) > 0 $,
\[
h^B_{top}( K_x , \mathcal{U} , \{F_n\} ) \ge f^{\{F_n\}}_{ \alpha , \mathcal{A} }(x) - \log M .
\]
In particular, if $ \mu $ is ergodic, then for $ \mu $-a.e. $ x \in X $ and any such $K_x$,
\[
h^B_{top} ( K_x , \mathcal{U} , \{F_n\} ) \geq h_{\mu} ( G , \alpha | \mathcal{A} ) - \log M.    
\]
\item[(2)] If $ \mu $ is ergodic, then for $\mu$-a.e. $ x \in X $ and any $ K_x \in \mathcal{B}_X $ with $ \mu_{x}(K_x) > 0 $, 
\[
h^B_{top} ( K_x , \{F_n\} ) \geq h_{\mu} (G,X | \mathcal{A} ) .    
\]
\end{enumerate}
\end{proposition}

\begin{proof}
(1) By Theorem \ref{thm-re-SMB}, for $ \mu $-a.e. $ x \in X $,
\[
\lim\limits_{n \to \infty} \frac{1}{|F_n|} I_{\mu} (\alpha_{F_n} | \mathcal{A} )(x) =  f^{\{F_n\}}_{ \alpha , \mathcal{A} }(x), 
\] 
and the function $f^{\{F_n\}}_{ \alpha , \mathcal{A} }(x) $ is $ \mathcal{A} $-measurable. Moreover, by (\ref{Exp-constant-fiber}), the following hold:
\begin{enumerate}
\item[(i)] For every $ B \in \alpha_{F_n} $ and for $ \mu $-a.e. $ x \in X $, 
\[
\mathbb{E}_{\mu} ( \mathbbm{1}_B | \mathcal{A} )(y) = \int_{X} \mathbbm{1}_B d \mu_y = \int_{X} \mathbbm{1}_B d \mu_x = \mathbb{E}_{\mu} ( \mathbbm{1}_B | \mathcal{A} )(x) \quad \text{for }\mu_x\text{-a.e. } y \in X.     
\] 
\item[(ii)] For $ \mu $-a.e. $ x \in X $,
\[
f^{\{F_n\}}_{ \alpha , \mathcal{A} }(y) = f^{\{F_n\}}_{ \alpha , \mathcal{A} }(x) \quad \text{for }\mu_x\text{-a.e. } y \in X.    
\] 
\end{enumerate}
Consequently, there exists a set $ X_0 \in \mathcal{B}_X $ with $ \mu(X_0)=1 $ such that for every $ x \in X_0 $, there is a measurable subset $ W_x \in \mathcal{B}_X $ with $ \mu_{x}(W_x) = 1 $ satisfying the following for all $ y \in W_x $:
\begin{enumerate}
\item[(iii)] For every $B \in \alpha_{F_n}$ and every $n \in \mathbb{N}$,
\[
\mathbb{E}_{\mu} ( \mathbbm{1}_B | \mathcal{A} )(y) = \mathbb{E}_{\mu} ( \mathbbm{1}_B | \mathcal{A} )(x) = \mu_{x}(B).    
\]
\item[(iv)] $\lim\limits_{n \to \infty} \frac{1}{|F_n|} I_{\mu} (\alpha_{F_n} | \mathcal{A} )(y) =  f^{\{F_n\}}_{ \alpha , \mathcal{A} }(y) =f^{\{F_n\}}_{ \alpha , \mathcal{A} }(x). $
\end{enumerate} 

Furthermore, for each $ x \in X_0 $ and $ y \in W_x $, we have
\begin{align}\label{equality-lowerable-prop}
\lim\limits_{n \to \infty} -\frac{1}{|F_n|} \log \mu_{x}( \alpha_{F_n}(y) ) &= \lim\limits_{n \to \infty} -\frac{1}{|F_n|} \log \mathbb{E}_{\mu} ( \mathbbm{1}_{ \alpha_{F_n}(y) } | \mathcal{A} )(y) \nonumber \\ &= \lim\limits_{n \to \infty} \frac{1}{|F_n|} I_{\mu} (\alpha_{F_n} | \mathcal{A} )(y) \nonumber \\ &= f^{ \{F_n\} }_{ \alpha , \mathcal{A} }(y).
\end{align}
		
Now fix $ x \in X_0 $ and let $ K_x \in \mathcal{B}_X $ be any subset with $ \mu_{x}(K_x) > 0 $. Then clearly $ \mu_{x}(W_x \cap K_x) = \mu_{x}(K_x) >0 $. 
For any $ \delta > 0 $ and $ n \in \mathbb{N} $, define 
\[
K_{x,l} (\delta) := \{ y \in W_x \cap K_x :  \mu_{x}( \alpha_{F_n}(y) ) \leq e^{ -|F_n|( f^{ \{F_n\}}_{ \alpha , \mathcal{A} }(x) - \delta ) } , \forall n \geq l \}.
\] 
By (\ref{equality-lowerable-prop}), we have $ K_x \cap W_x = \bigcup_{l \in \mathbb{N}} K_{x,l} (\delta) $. 
Since $ \mu_{x}(K_x \cap W_x) > 0 $, there exists $N \in \mathbb{N}$ such that $ \mu_{x}( K_{x,N}(\delta) ) > 0 $.

Note that for every $y \in K_{x,N}( \delta ) $, 
\[
\mu_{x}( \alpha_{F_n} (y) ) \leq e^{ -|F_n|( f^{ \{F_n\}}_{ \alpha , \mathcal{A} }(x) - \delta ) } \text{ , } \forall n \geq N .    
\] 
Therefore, for all $n \in \mathbb{N}$,
\[
\mu_{x}( \alpha_{F_n} (y) ) \leq c(y) \cdot e^{ -|F_n|( f^{ \{F_n\}}_{ \alpha , \mathcal{A} }(x) - \delta ) } \text{ , } \forall n \in \mathbb{N} ,    
\] 
where $ c(y) = \max \{ 1 , \sum_{i=1}^{N-1} e^{ |F_n|( f^{ \{F_i\}}_{ \alpha , \mathcal{A} }(x) - \delta ) } \} \in (0,\infty)$. 
By Lemma \ref{lem-nonuniform-massdis-principle}, it follows that
\[
h^B_{top}( K_{x,N}(\delta) , \mathcal{U} , \{F_n\} ) \geq f^{ \{F_n\} }_{ \alpha , \mathcal{A} }(x) - \delta -\log M .     
\] 
Hence, for every $\delta > 0$,
\[
h^B_{top}( K_x , \mathcal{U} , \{F_n\} ) \geq f^{ \{F_n\} }_{ \alpha , \mathcal{A} }(x) - \delta -\log M .    
\]
Since $\delta$ is arbitrary, we conclude that for all $x \in X_0$,
\[
h^B_{top}( K_x , \mathcal{U} , \{F_n\} ) \geq f^{ \{F_n\} }_{ \alpha , \mathcal{A} }(x) -\log M .    
\] 
		
In particular, if $ \mu $ is ergodic, then $ f^{ \{F_n\} }_{ \alpha , \mathcal{A} }(x) = h_{\mu} ( G , \alpha | \mathcal{A} ) $ for $ \mu $-a.e. $ x \in X $ and therefore
\[ 
h_{\mathrm{top}}^B (K_x , \mathcal{U}, \{F_n\}) \ge h_{\mu}(G,\alpha | \mathcal{A})- \log M .
\]

(2) Set $h := h_{\mu}(G, X | \mathcal{A})$. We aim to show that for $\mu$-a.e. $x \in X$ and any $K_x \in \mathcal{B}_X$ with $\mu_x(K_x) > 0$, one has
\[
h_{\mathrm{top}}^B (K_x, \{F_n\}) \ge h.
\]

By Theorem \ref{thm-re-BrinKatok}, there exists a measurable subset $X_1 \subseteq X$ with $\mu(X_1) = 1$ such that for every $x \in X_1$, there exists $X_2^x \in \mathcal{B}_X$ with $\mu_x(X_2^x) = 1$ and for every $y \in X_2^x$,
\[
\lim_{\varepsilon \to 0} \limsup_{n \to \infty} -\frac{\log \mu_x(B_{F_n}(y, \varepsilon))}{|F_n|} = h.
\]

Now, fix $x \in X_1$ and let $K_x \in \mathcal{B}_X$ with $\mu_x(K_x) > 0$. Set $a := \mu_x(K_x \cap X_2^x)$. Since $\mu_x(X_2^x) = 1$, we have $a = \mu_x(K_x) > 0$.

Let $\delta > 0$ be arbitrary (e.g., $\delta = 1$). By the pointwise limit above, for sufficiently small $\varepsilon > 0$, there exists a Borel subset $X_3 \subseteq X_2^x$ with $\mu_x(X_3) > 1 - \frac{a}{2}$ and $N \in \mathbb{N}$ such that for any $y \in X_3$ and any $n \ge N$,
\[
\mu_x(B_{F_n}(y, \varepsilon)) \le e^{-(h - \delta) \cdot |F_n|}.
\]
Note that $\mu_x(K_x \cap X_3) \ge \mu_x(K_x) - \mu_x(X \setminus X_3) > a - \frac{a}{2} = \frac{a}{2} > 0$.

We now estimate the Bowen topological entropy from below by constructing a suitable cover. Let $\{B_{F_{n_i}}(x_i, \frac{\varepsilon}{2})\}_i$ be a finite or countable family such that $x_i \in X$, $n_i \ge N$, and
\[
\bigcup_i B_{F_{n_i}}(x_i, \frac{\varepsilon}{2}) \supseteq K_x \cap X_3.
\]
We may assume without loss of generality that for each $i$, $B_{F_{n_i}}(x_i, \frac{\varepsilon}{2}) \cap K_x \cap X_3 \ne \emptyset$. Choose $y_i \in B_{F_{n_i}}(x_i, \frac{\varepsilon}{2}) \cap K_x \cap X_3$ for each $i$. Then by the triangle inequality, $B_{F_{n_i}}(x_i, \frac{\varepsilon}{2}) \subseteq B_{F_{n_i}}(y_i, \varepsilon)$.

Now, consider the Bowen outer measure $\mathcal{M}(K_x, N, \frac{\varepsilon}{2}, h-\delta, \{F_n\})$. We have
\begin{align*}
\mathcal{M}\left(K_x, N, \frac{\varepsilon}{2}, h-\delta, \{F_n\}\right) &\ge \sum_i \exp\left( -(h-\delta) \cdot |F_{n_i}| \right) \\
&\ge \sum_i \mu_x\left( B_{F_{n_i}}(y_i, \varepsilon) \right) \quad \text{(by the uniform estimate for $y_i \in X_3$)} \\
&\ge \sum_i \mu_x\left( B_{F_{n_i}}(x_i, \frac{\varepsilon}{2}) \right) \quad \text{(since $B_{F_{n_i}}(x_i, \frac{\varepsilon}{2}) \subseteq B_{F_{n_i}}(y_i, \varepsilon)$)} \\
&\ge \mu_x(K_x \cap X_3) \quad \text{(as the balls cover $K_x \cap X_3$)} \\
&> \frac{a}{2} > 0.
\end{align*}

Since this holds for all sufficiently small $\varepsilon$ and all $N$ large enough, we conclude that
\[
\mathcal{M}(K_x, h-\delta, \{F_n\}) = \lim_{\varepsilon \to 0} \lim_{N \to \infty} \mathcal{M}(K_x, N, \varepsilon, h-\delta, \{F_n\}) > 0.
\]
By the definition of Bowen topological entropy, this implies
\[
h_{\mathrm{top}}^B (K_x, \{F_n\}) \ge h - \delta.
\]
As $\delta > 0$ was arbitrary, we finally obtain
\[
h_{\mathrm{top}}^B (K_x, \{F_n\}) \ge h = h_{\mu}(G, X | \mathcal{A}),
\]
which completes the proof (2).   
\end{proof}

\begin{question}\label{rmk-construction-pf}
Can Proposition \ref{prop-lowerable-essential} (2) be proved utilizing the method in \cite{Dou2018Zhang}?
\end{question}

The following result is a direct application of Proposition \ref{prop-lowerable-essential}.
\begin{lemma}\label{lem-apply-lowerable-prop}
Let $(X,G)$ be a $G$-system, $\{F_n\} $ be a tempered F{\o}lner sequence of $G$ with $e_G \in F_1 \subseteq F_2 \subseteq \cdots $, and suppose it satisfies the growth condition $ \lim\limits_{n \to \infty } \frac{|F_n|}{\log n} = \infty $. Let $\mu \in \mathcal{M}_G^e(X)$ and $\mathcal{A} \subseteq \mathcal{B}_{\mu} $ a $G$-invariant sub-$\sigma$-algebra. If $ \mu = \int_{X} \mu_{x} d \mu (x) $ denotes the disintegration of $ \mu $ over $ \mathcal{A} $, then the following statements hold.
\begin{enumerate}
\item[(1)] If $\alpha \in \mathcal{P}_X$, then for $\mu$-a.e. $x \in X$, the following holds: for each fixed $x$ and any $\varepsilon \in (0,1)$, there exists a compact subset $K_x(\alpha,\varepsilon) \subseteq X$ such that $ \mu_{x}( K_x(\alpha,\varepsilon) ) \geq 1-\varepsilon $ and 
\[
h^B_{top}( K_x(\alpha,\varepsilon) , \alpha , \{F_n\} ) = h_{\mathrm{top}}( K_x(\alpha,\varepsilon) , \alpha , \{F_n\} ) = h_{\mu}(G,\alpha | \mathcal{A}) .    
\] 
\item[(2)] For $\mu$-a.e. $x \in X$, the following holds: for each fixed $x$ and any $\varepsilon \in (0,1)$, there exists a compact subset $K_x(\varepsilon) \subseteq X$ such that $\mu_{x}( K_x(\varepsilon) ) \geq 1-\varepsilon $ and 
\[
h^B_{top}( K_x(\varepsilon) , \{F_n\} ) = h_{\mathrm{top}}( K_x(\varepsilon) , \{F_n\} ) = h_{\mu}(G,X | \mathcal{A}) .    
\] 
\end{enumerate}
\end{lemma}

\begin{proof}
(1) Let $\alpha \in \mathcal{P}_X $. Because of $\mu \in \mathcal{M}_G^e (X) $, by (\ref{equality-lowerable-prop}), there exists $X_0 \in \mathcal{B}_X$ with $\mu (X_0) = 1 $ such that for each $x \in X_0 $, one can find $W_x \in \mathcal{B}_X$ with $\mu_{x} (W_x) = 1 $ and for each $ y \in W_x $, the following holds
\[
\lim\limits_{n\to \infty} -\frac{1}{| F_n |} \log \mu_{x} ( \alpha_{F_n} (y) ) = h_{\mu} ( G , \alpha | \mathcal{A} ) .    
\]  
By Proposition \ref{prop-lowerable-essential} (1), without loss of generality, for any $ x \in X_0 $ and $ K \in \mathcal{B}_X $ with $ \mu_{x} (K) > 0 $ (we can take $ K \cap X_0 $ as our consideration if necessary), one has 
\begin{align}\label{ineq-hBtop-lowerbound}
h^{B}_{top} ( K , \alpha , \{F_n\} ) \geq h_{\mu} ( G , \alpha | \mathcal{A} ).
\end{align}

Fixing $x \in X_0$, for each $y \in W_x$ and $m \in \mathbb{N}$, there exists $n(y, m) \in \mathbb{N}$ such that for all $n \ge n(y, m)$, we have
\[
-\frac{1}{|F_n|} \log \mu_{x} ( \alpha_{F_n}(y) ) \leq h_{\mu}(G,\alpha|\mathcal{A}) + \frac{1}{|F_m|}.
\]
Equivalently,
\begin{equation}\label{eq:entropy-bound-tail}
\mu_{x} ( \alpha_{F_n}(y) ) \geq e^{ -|F_n| \cdot ( h_{\mu}(G,\alpha|\mathcal{A}) + \frac{1}{|F_m|} ) } \quad \text{for all } n \ge n(y, m).
\end{equation}

We now define the function $c_m(y)$ for $\mu_x$-a.e. $y \in X$ by
\[
c_m(y) := \inf_{n \in \mathbb{N}} \frac{ \mu_{x} ( \alpha_{F_n}(y) ) }{ e^{ -|F_n| \cdot ( h_{\mu}(G,\alpha|\mathcal{A}) + \frac{1}{|F_m|} ) } }.
\]
Our goal is to show that $c_m(y) > 0$ for $y \in W_x$. To achieve this, we should analyze the following cases.

{\bf Case 1}: $n < n(y, m)$. 
    For these finitely many $n$, the quotient
    \[
    \frac{ \mu_{x} ( \alpha_{F_n}(y) ) }{ e^{ -|F_n| \cdot ( h_{\mu}(G,\alpha|\mathcal{A}) + \frac{1}{|F_m|} ) } }
    \]
    is positive. Let $M_1(y, m) > 0$ be the minimum value of this quotient over $n = 1, 2, \dots, n(y,m)-1$.

{\bf Case 2}: $n \ge n(y, m)$. 
    By inequality \eqref{eq:entropy-bound-tail}, the quotient is bounded below by $1$ for all $n$ in this range. Thus, we can take $M_2 = 1$.

Combining both cases, we conclude that for all $n \in \mathbb{N}$,
\[
\frac{ \mu_{x} ( \alpha_{F_n}(y) ) }{ e^{ -|F_n| \cdot ( h_{\mu}(G,\alpha|\mathcal{A}) + \frac{1}{|F_m|} ) } } \ge \min \{ M_1(y, m), 1 \} > 0.
\]
Taking the infimum over $n$, we obtain
\[
c_m(y) \ge \min \{ M_1(y, m), 1 \} > 0 \quad \text{for all } y \in W_x.
\]
Furthermore, $c_m(y)$ is $\mu_x$-measurable.

Now, for each $m, l \in \mathbb{N}$, define
\[
K(m,l) = \{ y \in W_x : c_m(y) \geq 1 / l \}.
\]
The set $K(m,l)$ is $\mu_x$-measurable. Since $c_m(y) > 0$ on $W_x$, we have $\bigcup_{l \in \mathbb{N}} K(m,l) = W_x$, and hence $\lim\limits_{l \to \infty} \mu_x(K(m,l)) = \mu_x(W_x) = 1$.

Finally, the uniform lower bound $c_m(y) \ge 1/l$ for $y \in K(m,l)$ allows us to apply Lemma \ref{lem-uniform-massdis-principle}, yielding 
\[
h_{\mathrm{top}}( K(m,l), \alpha, \{F_n\} ) \leq h_{\mu} ( G, \alpha | \mathcal{A} ) + \frac{1}{|F_m|}.
\]

Note that $ \lim_{l \to \infty} \mu_{x}( K(m,l) ) = 1 $. Hence, for each $ \varepsilon \in (0,1) $ and each $ m \in \mathbb{N} $, there exists some $ L_m \in \mathbb{N} $ such that $ \mu_x(K(m, L_m)) > 1 - \frac{\varepsilon}{2^{m+1}} $. By the inner regularity of the measure $ \mu_x $, we can find a compact subset $ B(m,\varepsilon) \subseteq K(m, L_m) $ satisfying
\[
\mu_{x} ( X \setminus B(m,\varepsilon) ) < \frac{\varepsilon}{2^m} \quad \text{and} \quad h_{\mathrm{top}}( B(m,\varepsilon) , \alpha , \{F_n\} ) \leq h_{\mu}( G , \alpha | \mathcal{A} ) + \frac{1}{|F_m|}.
\]
The entropy estimate for $ B(m,\varepsilon) $ follows from the fact that $ B(m,\varepsilon) \subseteq K(m, L_m) $ and Lemma \ref{lem-uniform-massdis-principle}.

Now, define
\[
K_x (\alpha,\varepsilon) = \bigcap_{m \in \mathbb{N} } B( m , \varepsilon ).
\]
Then $ K_x (\alpha , \varepsilon) $ is a compact subset of $ X $. We now verify its properties.

First, it has large measure, i.e.,
\begin{align*}
\mu_{x}( K_x (\alpha,\varepsilon) ) &= 1 - \mu_{x} \left( X \setminus \bigcap_{m \in \mathbb{N} } B(m,\varepsilon) \right) 
= 1 - \mu_{x} \left( \bigcup_{m \in \mathbb{N} } (X \setminus B(m,\varepsilon)) \right) \\
&\geq 1 - \sum_{m=1}^{\infty} \mu_{x} ( X \setminus B(m,\varepsilon) ) > 1 - \sum_{m=1}^{\infty} \frac{\varepsilon}{2^m} = 1 - \varepsilon > 0.
\end{align*}

Second, its topological entropy is bounded above by the conditional measure-theoretic entropy. Since $ K_x(\alpha,\varepsilon) \subseteq B(m,\varepsilon) $ for every $ m \in \mathbb{N} $, we have
\begin{align*}
h_{\mathrm{top}} ( K_x (\alpha,\varepsilon) , \alpha , \{F_n\} ) &\leq \inf_{m \in \mathbb{N} } h_{\mathrm{top}} ( B (m,\varepsilon) , \alpha , \{F_n\} ) \\
&\leq \inf_{m \in \mathbb{N} } \left( h_{\mu}(G,\alpha|\mathcal{A}) + \frac{1}{|F_m|} \right) \\
&= h_{\mu}(G,\alpha|\mathcal{A}). \quad \text{(since $ \lim_{m\to\infty} |F_m| = \infty $)}
\end{align*}

Finally, to establish the equality, we combine two results. On one hand, Lemma \ref{lem3.1} gives the general inequality $ h^{B}_{top}(\cdot,\{F_n\}) \leq h_{\mathrm{top}}(\cdot,\{F_n\}) $. On the other hand, the previously established inequality (\ref{ineq-hBtop-lowerbound}) provides 
\[
h_{\mu}(G, \alpha | \mathcal{A}) \leq h^{B}_{top}( K_x(\alpha,\varepsilon), \alpha, \{F_n\} ).
\]
Therefore,
\[
h_{\mu}(G, \alpha | \mathcal{A}) \leq h^{B}_{top}( K_x(\alpha,\varepsilon), \alpha, \{F_n\} ) \leq h_{\mathrm{top}}( K_x(\alpha,\varepsilon), \alpha, \{F_n\} ) \leq h_{\mu}(G, \alpha | \mathcal{A}),
\]
and we have
\[
h_{\mathrm{top}}^B ( K_x (\alpha,\varepsilon) , \alpha , \{F_n\} ) = h_{\mathrm{top}} ( K_x (\alpha,\varepsilon) , \alpha , \{F_n\} ) = h_{\mu}(G,\alpha|\mathcal{A}).
\]
This completes the proof of (1).

(2) Let $\{ \mathcal{U}_m \}_{m \in \mathbb{N}} \subseteq \mathcal{C}_{X}^o$ be a sequence of open covers with $\lim\limits_{m \to \infty} \| \mathcal{U}_m \| = 0$. For each $m \in \mathbb{N}$, let $\alpha_m \in \mathcal{P}_X$ be a finite Borel partition such that $\alpha_m \succeq \mathcal{U}_m$.

By the result of part (1), for each $m \in \mathbb{N}$, there exists a measurable subset $X_0^m \subseteq X$ with $\mu(X_0^m) = 1$ such that for every $x \in X_0^m$, every $\varepsilon \in (0,1)$, there exists a compact subset $K_x(m, \varepsilon)$ satisfying $\mu_x(K_x(m, \varepsilon)) \geq 1 - \frac{\varepsilon}{2^m}$ and
\[
h_{\mathrm{top}}^B ( K_x(m, \varepsilon), \alpha_m, \{F_n\} ) = h_{\mathrm{top}} ( K_x(m, \varepsilon), \alpha_m, \{F_n\} ) = h_{\mu}(G, \alpha_m | \mathcal{A}).
\]
Define $X_0 = \bigcap_{m=1}^{\infty} X_0^m$. Then $\mu(X_0) = 1$.

Now, fix $x \in X_0$ and $\varepsilon \in (0,1)$. Define
\[
K_x(\varepsilon) = \bigcap_{m \in \mathbb{N}} K_x(m, \varepsilon).
\]
Then $K_x(\varepsilon)$ is a compact subset of $X$. We verify its properties.

First, it has large measure, i.e.,
\begin{align*}
\mu_{x}( K_x(\varepsilon) ) &= 1 - \mu_{x} \left( X \setminus \bigcap_{m \in \mathbb{N}} K_x (m,\varepsilon) \right) 
= 1 - \mu_{x} \left( \bigcup_{m \in \mathbb{N}} (X \setminus K_x (m,\varepsilon)) \right) \\
&\geq 1 - \sum_{m=1}^{\infty} \mu_{x} ( X \setminus K_x (m,\varepsilon) ) 
\geq 1 - \sum_{m=1}^{\infty} \frac{\varepsilon}{2^m} = 1 - \varepsilon > 0.
\end{align*}

Second, we estimate its topological entropy. Recall that for any compact set $K$, $h_{\mathrm{top}}(K, \{F_n\}) = \sup_{\mathcal{U} \in \mathcal{C}_X^o} h_{\mathrm{top}}(K, \mathcal{U}, \{F_n\})$. Since $\{\mathcal{U}_m\}_{m \in \mathbb{N}}$ is a refining sequence of open covers with mesh going to zero, one has
\[
h_{\mathrm{top}} (K_x (\varepsilon), \{F_n\} ) = \sup_{m \in \mathbb{N}} h_{\mathrm{top}} ( K_x (\varepsilon), \mathcal{U}_m, \{F_n\} ).
\]
For each $m \in \mathbb{N}$, since $K_x(\varepsilon) \subseteq K_x(m, \varepsilon)$ and $\alpha_m \succeq \mathcal{U}_m$, we have
\begin{align*}
h_{\mathrm{top}} ( K_x (\varepsilon), \mathcal{U}_m, \{F_n\} ) &\leq h_{\mathrm{top}} ( K_x (m, \varepsilon), \mathcal{U}_m, \{F_n\} ) \\
&\leq h_{\mathrm{top}} ( K_x (m, \varepsilon), \alpha_m, \{F_n\} ) \\
&= h_{\mu} (G, \alpha_m | \mathcal{A}).
\end{align*}
Taking the supremum over $m$ on both sides, we obtain
\[
h_{\mathrm{top}} (K_x (\varepsilon), \{F_n\} ) \leq \sup_{m \in \mathbb{N}} h_{\mu} (G, \alpha_m | \mathcal{A}) \leq h_{\mu} (G, X | \mathcal{A}).
\]

Finally, to establish the equality, we combine the general upper bound from Lemma \ref{lem3.1}, $h^{B}_{top}(\cdot,\{F_n\}) \leq h_{\mathrm{top}}(\cdot,\{F_n\})$, with the previously established lower bound from Proposition \ref{prop-lowerable-essential} (2),
\[
h_{\mu}(G, X | \mathcal{A}) \leq h^{B}_{top}( K_x(\varepsilon), \{F_n\} ).
\]
Therefore,
\[
h_{\mu}(G, X | \mathcal{A}) \leq h^{B}_{top}( K_x(\varepsilon), \{F_n\} ) \leq h_{\mathrm{top}}( K_x(\varepsilon), \{F_n\} ) \leq h_{\mu}(G, X | \mathcal{A}),
\]
and we have
\[
h_{\mathrm{top}}^B (K_x (\varepsilon), \{F_n\} ) = h_{\mathrm{top}} (K_x (\varepsilon), \{F_n\} ) = h_{\mu} (G, X | \mathcal{A}).
\]
This completes the proof of (2).
\end{proof}
	
With the above preparations, we can now state and prove the main result of this section.
\begin{proof}[Proof of Theorem~\ref{thm-main-finite}]
If $h = h_{\mathrm{top}} (X,G)$, the statement holds trivially by taking $K_h = X$, as confirmed by Proposition \ref{prop-topoent-Bowentopoent}. 
If $h = 0$, it suffices to take $K_h = \{x\}$ for any $x \in X$, since the topological entropy of a single point is zero.

Now, assume $ 0 < h < h_{\mathrm{top}}(X,G) $. By the variational principle for amenable group actions, there exists an ergodic measure $ \mu \in \mathcal{M}_G^e(X) $ such that
\[
h < h_{\mu} (G,X) \leq h_{\mathrm{top}} (X,G) < \infty .   
\] 
By Proposition \ref{prop-exist-factor}, there exists a $G$-invariant sub-$ \sigma $-algebra $\mathcal{A} \subseteq \mathcal{B}_{\mu}$ such that 
\[
h_{\mu} (G,X|\mathcal{A}) = h .    
\] 
Now, applying Lemma \ref{lem-apply-lowerable-prop} (2) to the measure $\mu$ and the sub-$\sigma$-algebra $\mathcal{A}$, we obtain a compact subset $K_h\subseteq X$ satisfying
\[ 
h_{\mathrm{top}}^B (K_h,\{F_n\})= h_{\mathrm{top}}(K_h,\{F_n\})= h_{\mu} (G,X|\mathcal{A}) = h.
\]
Finally, combining this result with Proposition \ref{prop-entropy-B<P<UC},  we conclude that
\[ 
h_{\mathrm{top}}^P (K_h,\{F_n\})=h .
\]
This completes the proof.
\end{proof}

\subsection{Hereditary uniform lowerability for expansive $G$-systems}\label{section-enpansivehul}
We now pass from the finite-entropy lowering theorem to a hereditary construction for expansive actions.  The key point is that expansiveness converts local separation estimates into prescribed subset entropy and permits the compact sets to be chosen with at most one limit point.

\begin{definition}
Let $(X,G)$ be a $G$-system with metric $d$. It is said to be \textit{expansive} if there exists an $r>0$ such that any two distinct points $x \neq y$ in $X$ are separated by the group action, meaning
\[
\sup_{g \in G}d(gx,gy) \ge r.    
\]  
In this case, $r$ is referred to as an \textit{expansive constant} for the system.
\end{definition}

Let $(X,G)$ be a $G$-system, $\{F_n\} $ a F{\o}lner sequence of $G$, and $E$ a closed subset of $X$. For every $\varepsilon > 0$ and $x \in E$, define 
\[ 
h_{\mathrm{spn}} (x,\varepsilon,E,\{F_n\}) = \inf \left\{h_{\mathrm{spn}}(d,\varepsilon,K) : K \text{ is a closed neighborhood of } x \text{ in } E \right\} .
\]
Let 
\[
h_{\mathrm{spn}}(x,E,\{F_n\}) = \lim_{\varepsilon \to 0^+} h_{\mathrm{spn}} (x,\varepsilon,E,\{F_n\}) ,    
\]
a quantity that depends only on the topology of $X$. The following properties were established in \cite[Remark~5.16]{Wang2026}.

\begin{proposition}\label{prop-general-hspn(x)}
Let $(X,G)$ be a $G$-system with metric $d$, $\{F_n\} $ a F{\o}lner sequence of $G$, and $E$ a closed subset of $X$. Then 
\begin{enumerate}
\item For every $\varepsilon> 0$, the function $x \mapsto h_{\mathrm{spn}}(x,\varepsilon,E,\{F_n\})$ is upper semi-continuous on $E$. Moreover,
\[
\sup_{x\in E} h_{\mathrm{spn}}(x,E,\{F_n\}) = h_{\mathrm{top}}(E,\{F_n\}) .
\]
\item If $E'$ is a countable closed subset of $E$ with a unique limit point $x$ in $E$, then  
\[
h_{\mathrm{top}}(E',\{F_n\}) \le h_{\mathrm{spn}}(x,E,\{F_n\}) .
\]
\item For every $x \in E$, there exists a countable closed subset $E_x \subseteq E$ containing $x$ such that $x$ is the unique limit point of $E_x$ in $E$ and 
\[
h_{\mathrm{top}}(E_x,\{F_n\})=h_{\mathrm{spn}}(x,E,\{F_n\}) .    
\]
\item There exists a countable closed subset $E' \subseteq E$ such that 
\[
h_{\mathrm{top}}(E',\{F_n\}) = h_{\mathrm{top}}(E,\{F_n\}) .   
\] 
Furthermore, $E'$ can be chosen so that the set of its limit points in $E$ has at most one limit point. In this case,  $E'$ has a unique limit point if and only if there exists some $x \in E$ satisfying
\[
h_{\mathrm{spn}}(x,E,\{F_n\}) = h_{\mathrm{top}}(E,\{F_n\}) .   
\]
\end{enumerate}
\end{proposition}

Note that this definition depends only on the topology of $X$, as any two compatible metrics on $X$ are uniformly equivalent.

\begin{lemma}\label{lem-expsys-lowofsubset}
Let $(X,G)$ be a $G$-system with metric $d$, $\{F_n\} $ a F{\o}lner sequence of $G$ satisfying
\[
e_G \in F_1 \subseteq F_2 \subseteq \cdots \quad \text{and} \quad \lim\limits_{n \to \infty} \frac{|F_n|}{\log n} = \infty,
\]
and let $K \subseteq X$ be a compact subset with $h_{\mathrm{top}}(K, \{F_n\}) > 0$.
Then for any $0 < h < h_{\mathrm{top}} (K, \{F_n\})$, there exists $ \delta_0 $ such that for every $ \delta \in (0,\delta_0] $, one can find a countable compact subset $K(h,\delta) \subseteq K$ with a unique limit point and with the following property
\[
h_{\mathrm{sep}}(d,\delta, K(h,\delta) ) = h .   
\] 
\end{lemma}

\begin{proof}
By Proposition \ref{prop-general-hspn(x)}, there exists a countable closed subset $K_0 \subseteq K$ such that $x_0 \in K_0$ is its unique limit point and $h_{\mathrm{top}}(K_0,\{F_n\}) >h$. Hence, there exists $ \delta_0 > 0 $ such that for every $0 < \delta \le \delta_0$, one has $ h_{\mathrm{sep}}(d,\delta,K_0) > h$.

Fix $0 < \delta \le \delta_0$. Define 
\begin{align*}
l_1 = \min \big\{ n \in \mathbb N :\, &\text{there exists } B_1 \subseteq K_0 \text{ such that }
\\ & B_1 \text{ is } (F_n , \delta)\text{-separated and } |B_1| = \left[e^{ h|F_n| } \right]+2  \big\}.
\end{align*}
Since $h_{\mathrm{sep}}(d,\delta,K_0) > h$, the number $l_1$ is finite. Choose $A_1 = D_1 \subseteq K_0$ to be an $(F_{l_1},\delta)$-separated set such that $|A_1| = \left[e^{ h|F_{l_1}| }\right]+ 1$ and $x_0 \notin A_1$.

Next, define
\begin{align*}
l_2 = \min \big\{ n \in \mathbb N : n > l_1 ,\, &\text{there exists } B_2 \subseteq \left( K_0 \cap B_{F_{l_1} }(x_0,\delta) \right) \setminus A_1  \text{ such that } \\ & B_2 \text{ is } (F_n , \delta)\text{-separated and } |B_2|=\left[e^{ h|F_n| }\right]-\left[ e^{h|F_{l_1}|} \right]+2  \big\}.
\end{align*}
Because $x_0$ is the unique limit point of $K_0$ in $K$, the set $K_0 \setminus \left( \left(K_0 \cap B_{F_{l_1}}(x_0,\delta) \right) \setminus A_1 \right) $ is finite. Therefore,
\[ 
h_{\mathrm{sep}} \left(d,\delta, \left(K_0 \cap B_{F_{l_1}}(x_0,\delta) \right) \setminus A_1 \right) = h_{\mathrm{sep}} (d, \delta , K_0) > h ,\]
which implies that $l_2 > l_1$ and $l_2$ is finite. Now, choose $D_2 \subseteq \left(K_0 \cap B_{F_{l_1}}(x_0, \delta) \right) \setminus A_1$ to be an $(F_{l_2}, \delta )$-separated set with $|D_2| = \left[e^{h |F_{l_2}|} \right] - \left[ e^{h |F_{l_1}| } \right] + 1 $ and $x_0 \notin D_2$. Define $A_2 = A_1 \cup D_2$. Then $x_0 \notin A_2$, $A_2 \subseteq K_0$ is $(F_{l_2},\delta)$-separated, and $|A_2| = \left[ e^{h |F_{l_2}|} \right] +2$.

Proceeding by induction, we obtain a strictly increasing sequence of natural numbers $\{l_i\}_{i \in \mathbb{N}}$ and an increasing sequence of subsets $\{A_i\}_{i \in \mathbb{N}}$ of $K_0$ such that for each $i \in \mathbb{N}$:
\begin{itemize}
\item $x_0 \notin A_i$ and $A_{i+1}\setminus A_i \subseteq K_0 \cap B_{ F_{l_i} }(x_0,\delta)$;
\item $A_i$ is $(F_{l_i},\delta)$-separated and $|A_i| = \left[ e^{h |F_{l_i}|} \right]+i$.
\end{itemize}

Now define $A = \left( \bigcup_{i \in \mathbb N} A_i \right) \cup \{x_0\} \subseteq K_0$. Then $A$ is a countable closed set with $x_0$ as its unique limit point in $X$. We claim that $h_{\mathrm{sep}}(d, \delta, A) = h$, and thus we may take $K(h, \delta) = A$. 

To verify the claim, fix $ l \in \mathbb N $ and choose $ n \in \mathbb{N}$ such that $l_n \le l \le l_{n+1}$. Let $\tilde{A} \subseteq A$ be an $(F_l,\delta)$-separated subset. Note that
\[
\tilde{A} \setminus A_n \subseteq \left( B_{ F_{l_n} }(x_0,\delta) \cap K_0 \right) \setminus A_n,    
\]
and this set remains $(F_l,\delta)$-separated. Therefore, 
\[ 
|\tilde{A} \setminus A_n| \le |A_{n+1} \setminus A_n| = \left[ e^{h |F_{l_{n+1}}|} \right] - \left[ e^{h |F_{l_n}|} \right] + 1. 
\]
It follows that
\begin{align*}
|\tilde{A}| &\le |\tilde{A} \setminus A_n| + |A_n| \\ & \le \left( \left[ e^{h |F_{l_{n+1}}|} \right] - \left[ e^{h |F_{l_n}|} \right] + 1 \right) + \left( \left[ e^{h |F_{l_n}|} \right] + n \right)\\ 
&= \left[ e^{h |F_{l_{n+1}}|} \right] + n + 1.
\end{align*}
Hence, $\mathrm{sep}(d , F_l , \delta , A) \le \left[ e^{h |F_{l_{n+1}}|} \right] + n + 1 $. 

On the other hand, since $A_n \subseteq A$ is $(F_{l_n}, \delta)$-separated, we have
\[
\mathrm{sep}(d, F_{l_n}, \delta, A) \ge |A_n| = \left[ e^{h |F_{l_n}|} \right] + n.
\]
Combining both estimates, we conclude that $h_{\mathrm{sep}}(d, \delta, A) = h$, which completes the proof.
\end{proof}

Since $G$ is amenable, we may fix a F{\o}lner sequence $\{F_n\} $ for the forthcoming theorem such that
\begin{align}\label{FS-condition}
e_G \in F_1 \subseteq F_2 \subseteq \cdots, \quad \lim\limits_{n \to \infty} \frac{|F_n|}{\log n} = \infty, \quad \lim_{n \to \infty} \frac{|F_{n+1}|}{|F_n|} = 1, \quad \text{and} \quad G = \bigcup_{n \in \mathbb{N}} F_n.   
\end{align}
The existence of such a special F{\o}lner sequence is guaranteed by Theorem 2.5 in \cite{Downarowicz2024OprochaZhang}, 
and it will be used in the next result without further explicit mention.

\begin{theorem}\label{thm-expansive-hul}
Let $(X,G)$ be an expansive $G$-system and $\{F_n\} $ a F{\o}lner sequence of $G$. Then $(X,G)$ is hereditarily uniformly lowerable along $\{F_n\} $. 
That is, for every non-empty compact subset $E \subseteq X$ and every $0 \le h \le h_{\mathrm{top}}(E, \{F_n\} )$, there exists a non-empty compact subset $K_h \subseteq K$ such that:
\begin{itemize}
    \item $h_{\mathrm{top}}(K_h, \{F_n\}) = h$,
    \item $K_h$ has at most one limit point.
\end{itemize}
\end{theorem}

\begin{proof}
Let $2r$ (with $r > 0$) be an expansive constant of the expansive $G$-system $(X, G)$.
Let $E \subseteq X$ be a non-empty compact subset and $0 \le h \le h_{\mathrm{top}}(E,\{F_n\})$. We consider three cases.

\textbf{Case (i)}: $h = 0$.  
For any $x \in E$, we have $h_{\mathrm{top}}(\{x\}, \{F_n\}) = 0$. Fix any $x \in E$ and take $K_h = \{x\}$.
    
\textbf{Case (ii)}: $h = h_{\mathrm{top}}(E,\{F_n\}) > 0$.  
By Proposition~\ref{prop-general-hspn(x)}, the function $x \mapsto h_{\mathrm{spn}} (x, r, E, \{F_n\})$ is upper semi-continuous on $E$ and 
\[
\sup_{x \in E} h_{\mathrm{spn}}(x, E, \{F_n\}) = h_{\mathrm{top}}(E,\{F_n\}).   
\]
Recall that
\begin{align*}
h_{\mathrm{spn}}(x, E,\{F_n\}) & = \lim_{\varepsilon \to 0^+ } h_{\mathrm{spn}}(x,\varepsilon,E,\{F_n\}) \\
& = \lim_{\varepsilon \to 0^+} \inf \left\{ h_{\mathrm{spn}}(d, \varepsilon, K) : K \text{ is a closed neighborhood of } x \text{ in } E \right\} \\
& = \lim_{\varepsilon \to 0^+} \inf \left\{ h_{\mathrm{sep}}(d, \varepsilon, K) : K \text{ is a closed neighborhood of } x \text{ in } E \right\} .
\end{align*}
The completion of this case relies on the following claim, which will be verified subsequently.

\begin{claim}\label{claim-thm-expansive-hul}
If $K \subseteq X$ is a compact subset, then $h_{\mathrm{top}}(K,\{F_n\}) = h_{\mathrm{sep}}(d, r, K)$.
\end{claim}

Assuming the claim for now, we obtain
\begin{align*}
h_{\mathrm{spn}}(x, E, \{F_n\}) &= \lim_{\varepsilon \to 0^+} \inf \left\{ h_{\mathrm{sep}}(d, r, K) : K \text{ is a compact neighborhood of } x \text{ in } E \right\} \\ &= h_{\mathrm{spn}} (x, r, E, \{F_n\}),
\end{align*}
and hence
\[ 
h_{\mathrm{top}}(E,\{F_n\}) = \max_{x \in E} h_{\mathrm{spn}}(x, r, E, \{F_n\}).  
\]

Choose $x_0 \in E$ such that $ h_{\mathrm{spn}}(x_0,r,E,\{F_n\}) = h_{\mathrm{top}}(E,\{F_n\}) >0$. Then by Proposition \ref{prop-general-hspn(x)}, there exists a countable closed subset $E_{x_0} \subseteq E$ with $x_0 \in E_{x_0}$ such that $x_0$ is the unique limit point of $E_{x_0}$ in $E$ and 
\begin{align*}
h_{\mathrm{top}}(E_{x_0}, \{F_n\}) &= h_{\mathrm{spn}}(x_0, E, \{F_n\}) \\
&=h_{\mathrm{spn}}(x_0, r, E, \{F_n\}) = h_{\mathrm{top}}(E, \{F_n\}) > 0.
\end{align*}
We therefore take $K_h = E_{x_0}$.

\textbf{Case (iii)}: $0 < h < h_{\mathrm{top}}(E, \{F_n\})$. By Lemma \ref{lem-expsys-lowofsubset}, there exists $\delta_0 > 0$ with $\delta_0 < r$ such that for every $0 < \delta < \delta_0$, there exists a countable compact subset $K(h,\delta) \subseteq E$ with a unique limit point satisfying $h_{\mathrm{sep}}(d, \delta, K(h,\delta)) = h$. In particular, $h_{\mathrm{sep}}(d, \delta_0 , K(h,\delta_0)) = h$. Since $2\delta_0$ is also an expansive constant, Claim~\ref{claim-thm-expansive-hul} implies that
\[
h_{\mathrm{top}}(K(h,\delta_0), \{F_n\}) = h,    
\]
and we take $K_h = K(h,\delta_0)$.
    
This establishes that for every non-empty compact subset $E \subseteq X$ and every $0 \le h \le h_{\mathrm{top}}(E, \{F_n\})$, there exists a non-empty compact subset $K_h \subseteq E$ with at most one limit point such that $h_{\mathrm{top}}(K_h, \{F_n\}) = h$.

It remains to validate Claim~\ref{claim-thm-expansive-hul}. Let $K \subseteq X$ be compact. We aim to show that $h_{\mathrm{top}}(K, \{F_n\}) = h_{\mathrm{sep}}(d, r, K)$. 
We first establish a key uniform continuity estimate: for every $\varepsilon > 0$, there exists $n(\varepsilon) \in \mathbb{N}$ such that for any $x, y \in X$ with $d_{F_{n(\varepsilon)}}(x, y) \le r$, one has $d(x, y) \le \varepsilon$.
Suppose this fails for some $\varepsilon > 0$. Then for every $n \in \mathbb{N}$, there exist $x_n, y_n \in X$ such that $d_{F_n}(x_n, y_n) \le r$ but $d(x_n, y_n) > \varepsilon$. By compactness, we may assume $(x_n, y_n) \to (x, y)$. Then $d(x, y) \ge \varepsilon$ and $d_{F_n}(x, y) \le r$ for all $n \in \mathbb{N}$, contradicting the fact that $2r$ is an expansive constant.

Now, fix $\varepsilon > 0$ and $m \in \mathbb{N}$. Let $E \subseteq K$ be an $(F_m, \varepsilon)$-separated set. We claim that if the F{\o}lner sequence satisfies the submultiplicativity condition $F_m F_{m'} \subseteq F_{m + m'}$, then $E$ is also $(F_{m + n(\varepsilon)}, r)$-separated.
Assume, for contradiction, that there exist distinct $x \neq y$ in $E$ with $d_{F_{m + n(\varepsilon)}}(x, y) \leq r$. 
For each $s \in F_m$, the inclusion $F_{n(\varepsilon)} s \subseteq F_{m + n(\varepsilon)}$ implies $d_{F_{n(\varepsilon)}}(sx, sy) \leq r$. By the uniform continuity estimate established above, this gives $d(sx, sy) \leq \varepsilon$. Hence, $d_{F_m}(x, y) \leq \varepsilon$, contradicting the $(F_m, \varepsilon)$-separation of $E$.
Therefore, under the submultiplicativity condition, we have
\begin{equation}
\mathrm{sep}(d, F_m, \varepsilon, K) \leq \mathrm{sep}(d, F_{m + n(\varepsilon)}, r, K).
\label{thm:Ineq-sep}
\end{equation}

For the original F{\o}lner sequence, we can select a subsequence $\{F_{m_L}\}_{L \in \mathbb{N}}$ satisfying $F_{m_L} F_{m_{L'}} \subseteq F_{m_{L + L'}}$ along which the entropy is attained
\[
h_{\mathrm{sep}}(d, \varepsilon, K) = \lim_{L \to \infty} \frac{1}{|F_{m_L}|} \log \mathrm{sep}(d, F_{m_L}, \varepsilon, K).
\]
Applying \eqref{thm:Ineq-sep} to this subsequence yields
\begin{align*}
h_{\mathrm{sep}}(d, \varepsilon, K)
&\leq \lim_{L \to \infty} \frac{1}{|F_{m_L}|} \log \mathrm{sep}(d, F_{m_L + n(\varepsilon)}, r, K) \\
&= \lim_{L \to \infty} \frac{1}{|F_{m_L + n(\varepsilon)}|} \log \mathrm{sep}(d, F_{m_L + n(\varepsilon)}, r, K) \quad 
\\
&\leq \limsup_{m \to \infty} \frac{1}{|F_m|} \log \mathrm{sep}(d, F_m, r, K) \\
&= h_{\mathrm{sep}}(d, r, K).
\end{align*}
Since $\varepsilon > 0$ was arbitrary and $h_{\mathrm{top}}(K, \{F_n\}) = \lim\limits_{\varepsilon \to 0^+} h_{\mathrm{sep}}(d, \varepsilon, K)$, we conclude that
\[
h_{\mathrm{top}}(K, \{F_n\}) \leq h_{\mathrm{sep}}(d, r, K).
\]
The reverse inequality is immediate from the definitions. This completes the proof.
\end{proof}


\section{Asymptotic $h$-expansiveness and D-hereditary lowerability}
\label{section-DHL-ahe}

In this section, we study hereditary lowering for Bowen topological entropy on arbitrary Souslin subsets. We first prove D-hereditary lowerability for finite-alphabet subshifts and then transfer this property to general asymptotically $h$-expansive actions by means of principal quasi-symbolic extensions. We also establish a quantitative lowering result for systems with positive topological tail entropy. Throughout this section, the F{\o}lner sequence $\{F_n\}$ is assumed to satisfy
\begin{equation}
	e_G\in F_1\subseteq F_2\subseteq\cdots.
	\label{eq:nested-F}
\end{equation}

\subsection{A factor inequality for Bowen entropy}
Let $\pi:(Y,G)\to(X,G)$ be a factor map.  We use the relative topological
entropy $h_{\mathrm{top}}(G,Y|\pi)$ introduced in
Subsection~\ref{subsection-tail-extension-tools}.
The fibre formula of Yan \cite{Yan2015} gives
\begin{equation}
	h_{\mathrm{top}}(G,Y|\pi)
	=\sup_{x\in X}h_{\mathrm{top}}(\pi^{-1}(x),\{F_n\}).
	\label{eq:fibre-formula-DHL}
\end{equation}

We shall use the following dimensional factor inequality.

\begin{lemma}\label{lem-factor-Bowen-DHL}
For every nonempty $E\subseteq Y$,
\begin{equation}
 h_{\mathrm{top}}^B(\pi E,\{F_n\})
 \le h_{\mathrm{top}}^B(E,\{F_n\})
 \le h_{\mathrm{top}}^B(\pi E,\{F_n\})
     +h_{\mathrm{top}}(G,Y|\pi).
 \label{eq:factor-inequality-DHL}
\end{equation}
The conclusion holds along every F{\o}lner sequence.
\end{lemma}

\begin{proof}
The first inequality follows from uniform continuity of $\pi$ and the
monotonicity of Bowen covers under factor maps.  The second is
\cite[Theorem~1.2]{HeZhang2026}, together with
\eqref{eq:fibre-formula-DHL}.
\end{proof}

\begin{corollary}\label{cor-zero-relative-DHL}
If $h_{\mathrm{top}}(G,Y|\pi)=0$, then
\[
 h_{\mathrm{top}}^B(E,\{F_n\})
 =h_{\mathrm{top}}^B(\pi E,\{F_n\})
\]
for every $E\subseteq Y$.  In particular, this holds when $\pi$ is a
principal extension.
\end{corollary}

\begin{proof}
Only the last assertion requires explanation.  For a principal extension the
relative measure-theoretic entropy is zero, by the Abramov--Rokhlin addition
formula for amenable actions \cite{WardZhang1992}.  The relative variational
principle \cite{Dooley2015Zhang,Yan2015} therefore gives
$h_{\mathrm{top}}(G,Y|\pi)=0$, and the conclusion follows from
Lemma~\ref{lem-factor-Bowen-DHL}.
\end{proof}

\begin{lemma}\label{lem-transfer-DHL}
Suppose $\pi:(Y,G)\to(X,G)$ satisfies
$h_{\mathrm{top}}(G,Y|\pi)=0$.  Then $(X,G)$ is D-hereditarily
lowerable along $\{F_n\}$ if and only if $(Y,G)$ is.
\end{lemma}

\begin{proof}
Assume first that $X$ is D-hereditarily lowerable and let $K\subseteq Y$ be
a nonempty Souslin set.  Then $\pi K$ is Souslin and
\[
 h_{\mathrm{top}}^B(K,\{F_n\})
 =h_{\mathrm{top}}^B(\pi K,\{F_n\})
\]
by Corollary~\ref{cor-zero-relative-DHL}.  Given an admissible $h$, choose
$L_h\subseteq\pi K$ with $h_{\mathrm{top}}^B(L_h,\{F_n\})=h$ and put
\[
 K_h=K\cap\pi^{-1}(L_h).
\]
Since $\pi K_h=L_h$, another application of
Corollary~\ref{cor-zero-relative-DHL} gives
$h_{\mathrm{top}}^B(K_h,\{F_n\})=h$.

Conversely, assume that $Y$ is D-hereditarily lowerable and let
$K\subseteq X$ be nonempty and Souslin.  Then $\pi^{-1}K$ is Souslin and has
the same Bowen entropy as $K$.  Choose $E_h\subseteq\pi^{-1}K$ with entropy
$h$ and set $K_h=\pi E_h$.  Corollary~\ref{cor-zero-relative-DHL} completes
the proof.
\end{proof}

\subsection{D-hereditary lowerability of symbolic systems}
Let $A$ be a finite alphabet with $|A|\ge2$ and let
$\Sigma\subseteq A^G$ be a $G$-subshift.  We use the right-shift convention
\[
 (gx)(u)=x(ug),\qquad g,u\in G.
\]
Let $\alpha$ be the clopen partition according to the coordinate at $e_G$,
and for a finite set $D\subseteq G$ put
\[
 \alpha_D=\bigvee_{d\in D}d^{-1}\alpha.
\]

\begin{lemma}\label{lem-generator-DHL}
For every $E\subseteq\Sigma$ and every finite $D\subseteq G$,
\[
 h_{\mathrm{top}}^B(E,\alpha_D,\{F_n\})
 =h_{\mathrm{top}}^B(E,\alpha,\{F_n\}).
\]
Consequently,
\begin{equation}
 h_{\mathrm{top}}^B(E,\{F_n\})
 =h_{\mathrm{top}}^B(E,\alpha,\{F_n\}).
 \label{eq:generator-DHL}
\end{equation}
\end{lemma}

\begin{proof}
Under the chosen shift convention,
$(\alpha_D)_{F_n}=\alpha_{DF_n}$.  Each $\alpha_{F_n}$-atom splits into at most
$
 |A|^{|DF_n\setminus F_n|}
$
$\alpha_{DF_n}$-atoms.  By the F{\o}lner property,
\[
 \frac{|DF_n\setminus F_n|}{|F_n|}\longrightarrow0.
\]
Given $\delta>0$, the splitting number is at most
$e^{\delta|F_n|}$ for all sufficiently large $n$.  Refining each member of a
Bowen--Carath\'eodory cover therefore gives
\[
 \mathcal M(E,\alpha_D,N,s,\{F_n\})
 \le
 \mathcal M(E,\alpha,N,s-\delta,\{F_n\})
\]
for all sufficiently large $N$.  Hence
\[
 h_{\mathrm{top}}^B(E,\alpha_D,\{F_n\})
 \le h_{\mathrm{top}}^B(E,\alpha,\{F_n\})+\delta.
\]
Letting $\delta\downarrow0$.  Every finite open cover of a subshift is refined
by a finite-coordinate cylinder partition $\alpha_D$, which proves
\eqref{eq:generator-DHL}.
\end{proof}

Put
\[
 S=\bigcup_{n\ge1}F_n
\]
and let $q:A^G\to A^S$ be coordinate restriction.  Repeated terms may be
deleted from the F{\o}lner sequence without changing the tail of the
Carath\'eodory construction, so we may assume that
$F_n\subsetneq F_{n+1}$.  For $u\ne v$ in $A^S$, let
\[
 m(u,v)=\max\{n\ge1:u|_{F_n}=v|_{F_n}\},
\]
where $m(u,v)=0$ if the two points already differ on $F_1$, and set
\[
 \rho(u,v)=
 \begin{cases}
 0,&u=v,\\
 1,&m(u,v)=0,\\
 e^{-|F_{m(u,v)}|},&m(u,v)\ge1.
 \end{cases}
\]

\begin{remark}\label{rem-layered-ultrametric-DHL}
The function $\rho$ is a well-defined metric compatible with the product topology on $A^S$.  
Indeed, if $u\ne v$, then they differ at some
$s\in S$, and hence at some finite level $F_N$, so $m(u,v)<\infty$.  If $u$
and $v$ agree on $F_a$ while $v$ and $w$ agree on $F_b$, then $u$ and $w$
agree on $F_{\min\{a,b\}}$; consequently
\[
 \rho(u,w)\le\max\{\rho(u,v),\rho(v,w)\}.
\]
Moreover, the closed $\rho$-ball centered at $u$ with radius $e^{-|F_n|}$ is exactly the $F_n$-cylinder containing $u$.
Since every finite subset of $S$ is contained in some $F_N$, these cylinders form a basis for the product
topology.  Notice that $\rho$ is a genuine metric on $A^S$, whereas its
pullback through $q$ need only be a pseudometric on $\Sigma$ when $S\ne G$.
\end{remark}

\begin{lemma}\label{lem-hausdorff-model-DHL}
For every $E\subseteq\Sigma$,
\begin{equation}
 h_{\mathrm{top}}^B(E,\alpha,\{F_n\})
 =\dim_H^{\rho}q(E).
 \label{eq:entropy-dimension-DHL}
\end{equation}
Consequently,
\begin{equation}
 h_{\mathrm{top}}^B(E,\{F_n\})
 =h_{\mathrm{top}}^B(E,\alpha,\{F_n\})
 =\dim_H^{\rho}q(E).
 \label{eq:full-entropy-dimension-DHL}
\end{equation}
\end{lemma}

\begin{proof}
For $u\in A^S$ and $n\ge1$, denote by
\[
 [u]_{F_n}=\{v\in A^S:v|_{F_n}=u|_{F_n}\}
\]
the $F_n$-cylinder containing $u$.  By the definition of the layered
ultrametric,
\begin{align*}
 \operatorname{diam}_{\rho}[u]_{F_n}=e^{-|F_n|}.
\end{align*}
Indeed, any two points in $[u]_{F_n}$ agree on $F_n$, so their
$\rho$-distance is at most $e^{-|F_n|}$.  Since
$F_n\subsetneq F_{n+1}$ and $|A|\ge2$, one can choose two configurations in
$[u]_{F_n}$ which differ at a coordinate of $F_{n+1}\setminus F_n$, showing
that the upper bound is attained.

We next observe that Hausdorff covers may be replaced by cylinder covers
without increasing their total cost.  Let $C\subseteq A^S$ contain at least
two points, and define
\[
 k(C)=\min\{m(u,v):u,v\in C,\ u\ne v\}.
\]
Then all points of $C$ agree on $F_{k(C)}$, whereas some pair fails to agree
on $F_{k(C)+1}$.  Hence $C$ is contained in an $F_{k(C)}$-cylinder and
\[
 \operatorname{diam}_{\rho}C=e^{-|F_{k(C)}|}.
\]
A singleton is contained in cylinders of arbitrarily large levels.  It
follows that every sufficiently fine Hausdorff cover of $q(E)$ can be
replaced by an $F_n$-cylinder cover whose diameters, and hence whose Hausdorff
covering costs, do not increase.

On the dynamical side, the atoms of
\[
 \alpha_{F_n}=\bigvee_{g\in F_n}g^{-1}\alpha
\]
are determined precisely by the $\alpha_{F_n}$-names.  Thus the
$\alpha_{F_n}$-atoms meeting $E$ correspond exactly, under the restriction
map $q$, to the $F_n$-cylinders meeting $q(E)$.  Moreover, 
\begin{align*}
 e^{-s|F_n|}
 =\bigl(e^{-|F_n|}\bigr)^s
 =\bigl(\operatorname{diam}_{\rho}[u]_{F_n}\bigr)^s.
\end{align*}
More explicitly, for $s\ge0$ and $N\in\mathbb N$, an admissible cover in the
definition of
$
 \mathcal M(E,\alpha,N,\\s,\{F_n\})
$
consists of atoms of $\alpha_{F_{n_i}}$ with $n_i\ge N$ and has total cost
$
 \sum_i e^{-s|F_{n_i}|}.
$

Under $q$, it becomes a cover of $q(E)$ by cylinders of diameters
$e^{-|F_{n_i}|}\le e^{-|F_N|}$ and has the same total cost in the definition
of the $s$-dimensional Hausdorff content.  Conversely, by the cylinder
replacement established above, every Hausdorff cover of $q(E)$ at a
sufficiently small scale may be replaced by such a cylinder cover, again
without increasing its cost.

Therefore the Carath\'eodory critical exponent defining
$h_{\mathrm{top}}^B(E,\alpha,\{F_n\})$ coincides with the Hausdorff critical
exponent defining $\dim_H^{\rho}q(E)$.  This proves
\eqref{eq:entropy-dimension-DHL}.  Finally,
Lemma~\ref{lem-generator-DHL} gives
\[
 h_{\mathrm{top}}^B(E,\{F_n\})
 =h_{\mathrm{top}}^B(E,\alpha,\{F_n\})
 =\dim_H^{\rho}q(E),
\]
which proves \eqref{eq:full-entropy-dimension-DHL}.
\end{proof}

\begin{theorem}\label{thm-subshift-DHL}
Every finite-alphabet $G$-subshift is D-hereditarily lowerable along every
F{\o}lner sequence satisfying \eqref{eq:nested-F}.
\end{theorem}

\begin{proof}
Let $K\subseteq\Sigma$ be a nonempty Souslin set.  Since $q$ is continuous,
$q(K)$ is Souslin, and by Lemmas~\ref{lem-generator-DHL} and
\ref{lem-hausdorff-model-DHL},
\[
 h_{\mathrm{top}}^B(K,\{F_n\})=\dim_H^{\rho}q(K).
\]
For $h=0$, take a singleton.  If
$0<h<\dim_H^{\rho}q(K)$, then
$\mathcal H_{\rho}^h(q(K))=\infty$.  Howroyd's theorem
\cite[Theorem~1]{Howroyd1995} gives a compact set $C_h\subseteq q(K)$ such
that
\[
 0<\mathcal H_{\rho}^h(C_h)<\infty.
\]
Hence $\dim_H^{\rho}C_h=h$.  Put
\[
 K_h=K\cap q^{-1}(C_h).
\]
Since $q(K_h)=C_h$,  \eqref{eq:entropy-dimension-DHL} gives
$h_{\mathrm{top}}^B(K_h,\{F_n\})=h$.  
When
$h=h_{\mathrm{top}}^B(K,\{F_n\})$, one takes $K_h=K$.  
The same argument handles
every finite $h$ when $h_{\mathrm{top}}^B(K,\{F_n\})=+\infty$, while the infinite endpoint is obtained simply by taking $K$ itself.
\end{proof}

\subsection{Asymptotic $h$-expansiveness implies D-hereditary lowerability}
The quasi-symbolic extension characterization in Lemma~\ref{lem-equicondition-asyhexpan} links the symbolic result to general asymptotically $h$-expansive action.

\begin{proof}[Proof of Theorem~\ref{thm-main-DHL}]
By Lemma~\ref{lem-equicondition-asyhexpan}, there exists a principal
quasi-symbolic extension
\[
 \pi:(Y,G)\longrightarrow(X,G),
 \qquad
 Y\subseteq\Sigma\times\mathbf T,
\]
where $\Sigma$ is a finite-alphabet $G$-subshift and $\mathbf T$ is a
zero-entropy tiling system.  Let
\[
 p=\operatorname{pr}_{\Sigma}|_Y,
 \qquad
 \Sigma_0=p(Y).
\]
Then $p:(Y,G)\to(\Sigma_0,G)$ is a factor map.

We first show that $p$ has zero relative topological entropy.  Let
\[
 q_{\mathbf T}=\operatorname{pr}_{\mathbf T}|_Y.
\]
For $z\in\Sigma_0$, set
\[
 \mathbf T_z=\{t\in\mathbf T:(z,t)\in Y\}.
\]
Then
\[
 p^{-1}(z)=\{z\}\times\mathbf T_z,
\]
and
\[
 q_z:p^{-1}(z)\longrightarrow\mathbf T_z,
 \qquad q_z(z,t)=t,
\]
is a homeomorphism.  Choose compatible metrics $d_\Sigma$ and
$d_{\mathbf T}$ and endow $Y$ with the restriction of the maximum product
metric on $\Sigma\times\mathbf T$.  For every finite $F\subseteq G$ and every
$(z,t),(z,s)\in p^{-1}(z)$,
\[
\begin{aligned}
 d_F^Y((z,t),(z,s))
 &=\max_{g\in F}\max\{d_\Sigma(gz,gz),d_{\mathbf T}(gt,gs)\}\\
 &=\max_{g\in F}d_{\mathbf T}(gt,gs)
 =d_F^{\mathbf T}(t,s).
\end{aligned}
\]
Thus $q_z$ preserves the finite-set Bowen metrics on the fibre.  
Since
$\mathbf T_z\subseteq\mathbf T$,
\[
h_{\mathrm{top}}(p^{-1}(z),\{F_n\})
=h_{\mathrm{top}}(\mathbf T_z,\{F_n\})
\le h_{\mathrm{top}}(\mathbf T,G)=0.
\]
The fibre formula \eqref{eq:fibre-formula-DHL} therefore gives
\[
 h_{\mathrm{top}}(G,Y| p)=0.
\]
Since $\Sigma_0$ is a finite-alphabet $G$-subshift, it is
D-hereditarily lowerable by Theorem~\ref{thm-subshift-DHL}.  Applying
Lemma~\ref{lem-transfer-DHL} to $p$, we conclude that $Y$ is
D-hereditarily lowerable.

Finally, $\pi$ is principal, so
$h_{\mathrm{top}}(G,Y|\pi)=0$ by
Corollary~\ref{cor-zero-relative-DHL}. 
Applying Lemma~\ref{lem-transfer-DHL} to $\pi$ gives the
D-hereditary lowerability of $(X,G)$.
\end{proof}

\begin{remark}\label{rem-nonnested-DHL}
The quasi-symbolic extension and the factor inequality are valid for
arbitrary F\o lner sequences, whereas nestedness is used in
Lemma~\ref{lem-hausdorff-model-DHL}.  For a general non-nested F\o lner
sequence, the argument reduces the problem to D-hereditary
lowerability of finite-alphabet full $G$-shifts along that sequence,
but the layered metric construction no longer applies directly.
Thus Theorem~\ref{thm-main-DHL} does not assert the unrestricted non-nested
version.
\end{remark}

\subsection{Tail entropy bounds for lowering Bowen entropy}
We now prove Theorem~\ref{thm-main-tailinterval}.  The argument uses a
principal zero-dimensional extension and finite symbolic factors whose
relative entropies decrease to the topological tail entropy.

\begin{lemma}\label{lem-zero-dimensional-DHL}
Every compact metrizable amenable $G$-action admits a principal
zero-dimensional extension
\[
 q:(Z,G)\longrightarrow(X,G).
\]
Moreover,
\[
 h_{\mathrm{top}}^*(Z,G)=h_{\mathrm{top}}^*(X,G).
\]
\end{lemma}

\begin{proof}
The existence of such an extension follows from
\cite[Theorem~3.2]{Huczek2021}, while invariance of topological tail entropy
under principal extensions is given by
\cite[Proposition~6.1]{DZ2022Tail}.
\end{proof}

Let $Z$ be zero-dimensional and choose finite clopen partitions
\[
 \mathcal P_1\prec\mathcal P_2\prec\cdots,
 \qquad
 \operatorname{mesh}(\mathcal P_j)\longrightarrow0.
\]
Let
\[
 \phi_j:(Z,G)\longrightarrow(Z_j,G)\subseteq A_j^G
\]
be the finite-alphabet symbolic factor generated by $\mathcal P_j$, and set
\[
 r_j=h_{\mathrm{top}}(G,Z|\phi_j).
\]

We use the tail variational principle of Downarowicz and Zhang
\cite[Theorem~2.1]{DZ2022Tail} in the following form.  If
$h_{\mathrm{top}}(Z,G)<+\infty$ and $\mu\in\mathcal M_G(Z)$, write
\[
 h_j(\mu)=h_\mu(G,\mathcal P_j),
 \qquad
 h(\mu)=h_\mu(Z,G),
 \qquad
 \theta_j(\mu)=h(\mu)-h_j(\mu).
\]
Then
\begin{equation}
 h_{\mathrm{top}}^*(Z,G)
 =\lim_{j\to\infty}\sup_{\mu\in\mathcal M_G(Z)}\theta_j(\mu).
 \label{eq:tail-variational-DHL}
\end{equation}

\begin{lemma}\label{lem-tail-approximation-DHL}
If $Z$ has finite topological entropy, then
\[
 r_j\downarrow h_{\mathrm{top}}^*(Z,G).
\]
\end{lemma}

\begin{proof}
By the relative variational principle in Proposition~\ref{prop-VP-yan},
\[
\begin{aligned}
 r_j
 &=\sup_{\mu\in\mathcal M_G(Z)}h_\mu(G,Z|\phi_j)\\
 &=\sup_{\mu\in\mathcal M_G(Z)}
   \bigl(h_\mu(Z,G)-h_{\phi_{j*}\mu}(Z_j,G)\bigr)\\
 &=\sup_{\mu\in\mathcal M_G(Z)}
   \bigl(h_\mu(Z,G)-h_\mu(G,\mathcal P_j)\bigr)
 =\sup_{\mu\in\mathcal M_G(Z)}\theta_j(\mu).
\end{aligned}
\]
The third equality follows because the zero-coordinate partition of $Z_j$ is
generating and its pullback under $\phi_j$ is $\mathcal P_j$.  Since the
partitions refine, $\{r_j\}$ is decreasing, and
\eqref{eq:tail-variational-DHL} gives the asserted limit.
\end{proof}

\begin{proof}[Proof of Theorem~\ref{thm-main-tailinterval}]
Take the principal zero-dimensional extension supplied by
Lemma~\ref{lem-zero-dimensional-DHL} and put
\[
 \widetilde K=q^{-1}K.
\]
Then $\widetilde K$ is Souslin and, by
Corollary~\ref{cor-zero-relative-DHL},
\[
 h_{\mathrm{top}}^B(\widetilde K,\{F_n\})=H,
 \qquad
 h_{\mathrm{top}}^*(Z,G)=a.
\]
For the symbolic factors above, Lemmas~\ref{lem-factor-Bowen-DHL} and
\ref{lem-tail-approximation-DHL} give
\begin{equation}
 h_{\mathrm{top}}^B(\phi_j\widetilde K,\{F_n\})\ge H-r_j,
 \qquad
 r_j\downarrow a.
 \label{eq:image-lower-bound-DHL}
\end{equation}

The case $h=0$ is attained by a singleton, so assume $h>0$.  If
$r_j\le a$ for some $j$, then
\[
 h_{\mathrm{top}}^B(\phi_j\widetilde K,\{F_n\})
 \ge H-r_j\ge H-a\ge h.
\]
By Theorem~\ref{thm-subshift-DHL}, choose
$L\subseteq\phi_j\widetilde K$ with
$h_{\mathrm{top}}^B(L,\{F_n\})=h$ and set
\[
 E=\widetilde K\cap\phi_j^{-1}L.
\]
Since $\phi_jE=L$, Lemma~\ref{lem-factor-Bowen-DHL} gives
\[
 h\le h_{\mathrm{top}}^B(E,\{F_n\})\le h+r_j\le h+a.
\]

It remains to consider the case $r_j>a$ for every $j$.  After discarding
finitely many terms, assume $r_j<a+h$ and define
\[
 \ell_j=h+a-r_j.
\]
Then $0<\ell_j<h$ and $\ell_j\uparrow h$.  Since $h+a\le H$,
\[
 \ell_j\le H-r_j
 \le h_{\mathrm{top}}^B(\phi_j\widetilde K,\{F_n\}).
\]
Choose $L_j\subseteq\phi_j\widetilde K$ with
$h_{\mathrm{top}}^B(L_j,\{F_n\})=\ell_j$ and set
\[
 E_j=\widetilde K\cap\phi_j^{-1}L_j.
\]
The factor inequality gives
\[
 \ell_j\le h_{\mathrm{top}}^B(E_j,\{F_n\})\le\ell_j+r_j=h+a.
\]
Let $E=\bigcup_{j\ge1}E_j$.  By countable stability of Bowen topological
entropy,
\[
 h_{\mathrm{top}}^B(E,\{F_n\})
 =\sup_{j\ge1}h_{\mathrm{top}}^B(E_j,\{F_n\}),
\]
and hence
\[
 h\le h_{\mathrm{top}}^B(E,\{F_n\})\le h+a.
\]
In either case set $K_h=q(E)\subseteq K$.  Since $q$ is principal,
Corollary~\ref{cor-zero-relative-DHL} yields
\[
 h_{\mathrm{top}}^B(K_h,\{F_n\})
 =h_{\mathrm{top}}^B(E,\{F_n\})\in[h,h+a].
\]
\end{proof}

\begin{remark}
	The adjustment $\ell_j=h+a-r_j$ is essential at the endpoint.  
	Although
	$r_j\downarrow a$, the Bowen entropies attained by the previously constructed
	sets need not form a closed set.  
	The modified targets ensure that all
	pullbacks remain below the common upper bound $h+a$, while the corresponding
	lower bounds converge to $h$.  
	Countable stability then yields the desired
	interval.
\end{remark}


\section{Hereditary uniform lowerability and asymptotic $h$-expansiveness}\label{section-asyhexpansive=hul}
In this section, we study the uniform hereditary lowering property for compact subsets with respect to the usual topological entropy.  The preceding section established D-hereditary lowerability for Bowen entropy under vanishing tail entropy.  Under the regularity assumptions on the F{\o}lner sequence imposed in Section~\ref{section-Finite-entropy-system}, we show that the stronger uniform hereditary lowering property is equivalent to asymptotic $h$-expansiveness.  The two implications are proved in Subsections~\ref{subsection-hul=>asyhex} and \ref{subsection-asyhex=>hul}, respectively.  We then examine the behavior of these lowering properties under principal extensions. 

\subsection{Asymptotic $h$-expansivity of hereditarily uniformly lowerable $G$-systems}\label{subsection-hul=>asyhex}
\begin{lemma}\label{lem-existzeroentropy-subset}
Let $(X,G)$ be a $G$-system with metric $d$, and let  $\{F_n\} $ be a F{\o}lner sequence of $G$ satisfying
\[
e_G \in F_1 \subseteq F_2 \subseteq \cdots \quad \text{and} \quad \lim\limits_{n \to \infty} \frac{|F_n|}{\log n} = \infty.
\]
If $K \subseteq X $ is a compact subset satisfying $h_{\mathrm{top}} (K,\{F_n\}) > 0$, then there exists an infinite countable compact subset $K_0 \subseteq K $ such that 
\[ 
h_{\mathrm{top}} (K_0,\{F_n\}) = 0 . 
\]
\end{lemma}

\begin{proof}
By Proposition~\ref{prop-general-hspn(x)}, there exists an infinite countable compact subset $\tilde{K} = \{x_0, x_1,\allowbreak x_2, \dots\} \subseteq K$ such that
\[
h:= h_{\mathrm{top}}(\tilde{K}, \{F_n \}) > 0 \quad \text{and} \quad \lim_{n \to \infty} x_n = x_0. 
\]

We now construct a nested sequence of subsets by induction. By Lemma~\ref{lem-expsys-lowofsubset}, for sufficiently small $\varepsilon_1 > 0$, there exists an infinite countable compact subset $K_1 \subseteq \tilde{K}$ such that
\[
h_{\mathrm{sep}}(d, \varepsilon_1, K_1) = \frac{h}{|F_1|}.
\]

Suppose $K_n$ has been constructed. Then for sufficiently small $\varepsilon_{n+1} > 0$ with $\varepsilon_{n+1} < \varepsilon_n$, applying Lemma~\ref{lem-expsys-lowofsubset} again yields an infinite countable compact subset $K_{n+1} \subseteq K_n$ such that
\[
h_{\mathrm{sep}}(d, \varepsilon_{n+1}, K_{n+1}) = \frac{h}{|F_n|}.
\]
We may choose $\{\varepsilon_n\}_{n \in \mathbb{N}}$ such that $\lim\limits_{n \to \infty} \varepsilon_n = 0$.

Now define $K_0 = \{x_0, \tilde{x}_1, \tilde{x}_2, \dots\}$ by selecting $\tilde{x}_n \in K_n \setminus K_{n+1}$ for each $n \in \mathbb{N}$. Then $K_0$ is an infinite countable compact subset. Moreover, for every $n \ge 1$,
\[
h_{\mathrm{sep}}(d, \varepsilon_n, K_0) \le h_{\mathrm{sep}}(d, \varepsilon_n, K_n) = \frac{h}{|F_n|}.
\]
Therefore,
\[
h_{\mathrm{top}}(K_0, \{F_n\}) = \lim_{n \to \infty} h_{\mathrm{sep}}(d, \varepsilon_n, K_0) = 0,
\]
which completes the proof.
\end{proof}

We first prove one direction of the equivalence in Theorem~\ref{thm-main-HUL}.

\begin{theorem}\label{thm-hul->asyhex}
Let $(X,G)$ be a $G$-system with metric $d$, and let $\{F_n\} $ be a F{\o}lner sequence satisfying
\[
e_G \in F_1 \subseteq F_2 \subseteq \cdots \quad \text{and} \quad \lim\limits_{n \to \infty} \frac{|F_n|}{\log n} = \infty.
\]
If $(X,G)$ is hereditarily uniformly lowerable along $\{F_n\} $, then it is asymptotically $h$-expansive.
\end{theorem}

\begin{proof}
By Lemma~\ref{prop-character-tailentropy}, it suffices to prove that
\begin{align*}
\lim_{\varepsilon \to 0^+} \sup_{x \in X} h_{\mathrm{top}}(\Gamma_{\varepsilon, \{F_n\}}(x), \{F_n\})= \lim_{\varepsilon \to 0^+} \sup_{x \in X} h_{\mathrm{top}}(\Gamma_{\varepsilon}(x), \{F_n\}) = 0.   
\end{align*}
Suppose, for contradiction, that
\[
h := \lim_{\varepsilon \to 0^+} \sup_{x \in X} h_{\mathrm{top}}(\Gamma_{\varepsilon, \{F_n\}}(x), \{F_n\}) = \lim_{\varepsilon \to 0^+} \sup_{x \in X} h_{\mathrm{top}}(\Gamma_{\varepsilon}(x), \{F_n\}) > 0.
\]
Then there exists a sequence $\{\varepsilon_i\}_{i \in \mathbb{N}}$ with $\lim_{i \to \infty} \varepsilon_i = 0$ such that
\[
h = \lim_{i \to \infty} h_i, \quad \text{where } h_i := \sup_{x \in X} h_{\mathrm{top}}(\Gamma_{\varepsilon_i}(x), \{F_n\}).
\]
For any $\sigma > 0$, there exists $x_i \in X$ such that
\[
h_{\mathrm{top}}(\Gamma_{\varepsilon_i}(x_i), \{F_n\}) \ge h_i - \sigma.
\]
Without loss of generality, assume $x_i \to x_0$ as $i \to \infty$. Then
\[
\lim_{i \to \infty} h_{\mathrm{top}}(\Gamma_{\varepsilon_i, \{F_n\}}(x_i), \{F_n\})= \lim_{\varepsilon \to 0^+} \sup_{x \in X} h_{\mathrm{top}}(\Gamma_{\varepsilon}(x), \{F_n\}) = h > 0.
\]

If the set $\{i \in \mathbb{N} : x_0 \notin \Gamma_{\varepsilon_i}(x_i)\}$ is infinite, we may assume, by passing to a subsequence, that $x_0 \notin \Gamma_{\varepsilon_i}(x_i)$ for all $i \ge 1$. Since $(X,G)$ is hereditarily uniformly lowerable along $\{F_n\} $, for each $i$ there exists an infinite countable compact subset $X_i \subseteq \Gamma_{\varepsilon_i}(x_i)$ with unique limit point $y_i$ in $X$ such that
\[
a_i = h_{\mathrm{top}}(X_i, \{F_n\}) < h \quad \text{and} \quad \lim_{i \to \infty} a_i = h.
\]
Then $y_i \to x_0$ as $i \to \infty$. By Lemma~\ref{lem-existzeroentropy-subset}, we may assume
\[
h_{\mathrm{top}}(\{x_0, y_1, y_2, \dots\}, \{F_n\}) = 0.
\]

We now state a key technical claim.

\begin{claim}\label{claim-thmasyhex-technical}
Let $\{B_i\}_{i \in \mathbb{N}}$ be a sequence of compact subsets of $X$ with $\lim_{i \to \infty} B_i = \{x_0\}$ in the Hausdorff metric $d_H$, and suppose
\begin{equation}\label{claim-ContractionCondition}
\inf_{N \in \mathbb{N}} \lim_{i \to \infty} \sup_{n \ge N} \left[ \max_{g \in F_n} \mathrm{diam}(g B_i) \right] = 0.
\end{equation}
For any choice of $x_i \in B_i$, we have
\[
h_{\mathrm{top}}\left( \bigcup_{i \in \mathbb{N}} B_i \cup \{x_0\}, \{F_n\} \right) = \max \left\{ \sup_{i \in \mathbb{N}} h_{\mathrm{top}}(B_i, \{F_n\}),\ h_{\mathrm{top}}(\{x_0, x_1, x_2, \dots\}, \{F_n\}) \right\}.
\]
\end{claim}

\begin{proof}[Proof of Claim~\ref{claim-thmasyhex-technical}]
Let $\varepsilon > 0$. For each $n \in \mathbb{N}$, let $E_n$ be an $(F_n, \varepsilon)$-separated set for $\bigcup_{i \in \mathbb{N}} B_i \cup \{x_0\}$ with cardinality $|E_n| = \mathrm{sep}(d, F_n, \varepsilon, \bigcup_{i \in \mathbb{N}} B_i \cup \{x_0\})$.

By condition~\eqref{claim-ContractionCondition}, there exist $N_\varepsilon, I_\varepsilon \in \mathbb{N}$ such that for all $n \ge N_\varepsilon$ and $i \ge I_\varepsilon$,
\[
\max_{g \in F_n} \mathrm{diam}(g B_i) < \frac{\varepsilon}{10}.
\]
For fixed $n \ge N_\varepsilon$ and $i \ge I_\varepsilon$, we have $|E_n \cap B_i| \le 1$. Indeed, if $y \ne z \in E_n \cap B_i$, then
\[
d_{F_n}(y, z) = \max_{g \in F_n} d(gy, gz) < \frac{\varepsilon}{10},
\]
contradicting the $\varepsilon$-separation of $E_n$.

Let $E_n \cap \left( \bigcup_{i \ge I_\varepsilon} B_i \right) = \{z_1, \dots, z_\kappa\}$. For each $z_k$, there exists a unique $i_k \ge I_\varepsilon$ such that $z_k \in B_{i_k}$. Define $E'_n = \{x_{i_1}, \dots, x_{i_\kappa}\}$. For distinct $k_1, k_2 \in \{1, \dots, \kappa\}$,
\begin{align*}
d_{F_n}(x_{i_{k_1}}, x_{i_{k_2}}) &\ge d_{F_n}(z_{k_1}, z_{k_2}) - d_{F_n}(z_{k_1}, x_{i_{k_1}}) - d_{F_n}(z_{k_2}, x_{i_{k_2}}) \\
&> \varepsilon - \frac{\varepsilon}{5} = \frac{4\varepsilon}{5} > \frac{\varepsilon}{2},
\end{align*}
so $E'_n$ is $(F_n, \frac{\varepsilon}{2})$-separated in $\{x_0, x_1, x_2, \dots\}$.
Since $E_n \cap \left( \bigcup_{1 \le i < I_\varepsilon} B_i \right)$ is $(F_n, \varepsilon)$-separated in $\bigcup_{1 \le i < I_\varepsilon} B_i$, for $n \ge N_\varepsilon$ we have
\begin{align*}
\mathrm{sep}(d, F_n, \varepsilon, \textstyle\bigcup_{i \ge 1} B_i \cup \{x_0\}) 
&= |E_n| \\
&\le \left| E_n \cap \left( \textstyle\bigcup_{1 \le i < I_\varepsilon} B_i \right) \right| + \left| E_n \cap \left( \textstyle\bigcup_{i \ge I_\varepsilon} B_i \right) \right| \\
&\le \mathrm{sep}(d, F_n, \varepsilon, \textstyle\bigcup_{1 \le i < I_\varepsilon} B_i) + \kappa \\
&= \mathrm{sep}(d, F_n, \varepsilon, \textstyle\bigcup_{1 \le i < I_\varepsilon} B_i) + |E'_n| \\
&\le \mathrm{sep}(d, F_n, \varepsilon, \textstyle\bigcup_{1 \le i < I_\varepsilon} B_i) + \mathrm{sep}(d, F_n, \tfrac{\varepsilon}{2}, \{x_0, x_1, x_2, \dots\}).
\end{align*}
Therefore,
\begin{align*}
&h_{\mathrm{sep}}(d, \varepsilon, \textstyle\bigcup_{i \ge 1} B_i \cup \{x_0\}) \\
&\quad \le \max \left\{ h_{\mathrm{sep}}(d, \varepsilon, \textstyle\bigcup_{1 \le i < I_\varepsilon} B_i),\ h_{\mathrm{sep}}(d, \tfrac{\varepsilon}{2}, \{x_0, x_1, x_2, \dots\}) \right\} \\
&\quad \le \max \left\{ \max_{1 \le i < I_\varepsilon} h_{\mathrm{top}}(B_i, \{F_n\}),\ h_{\mathrm{top}}(\{x_0, x_1, x_2, \dots\}, \{F_n\}) \right\}.
\end{align*}
Since $\varepsilon > 0$ was arbitrary, the claim follows.
\end{proof}

The sequence $\{X_i\}_{i \in \mathbb{N}}$ satisfies condition~\eqref{claim-ContractionCondition}, so by Claim~\ref{claim-thmasyhex-technical},
\[
h_{\mathrm{top}}\left( \textstyle\bigcup_{i \in \mathbb{N}} X_i \cup \{x_0\}, \{F_n\} \right) = \max \left\{ \sup_{i \in \mathbb{N}} h_{\mathrm{top}}(X_i, \{F_n\}),\ h_{\mathrm{top}}(\{x_0, y_1, y_2, \dots\}, \{F_n\}) \right\} = h.
\]
Since $(X,G)$ is hereditarily uniformly lowerable along $\{F_n\} $, there exists an infinite countable compact subset $K \subseteq \bigcup_{i \in \mathbb{N}} X_i \cup \{x_0\}$ with unique limit point $z$ in $X$ such that
\[
h_{\mathrm{top}}(K, \{F_n\}) = h > 0.
\]

We now consider two cases.

\textbf{Case 1}: $z = x_0$. 
Then $K_i := K \cap X_i$ is finite for each $i \ge 1$. By Claim~\ref{claim-thmasyhex-technical},
\begin{align*}
h_{\mathrm{top}}(K, \{F_n\}) &\le h_{\mathrm{top}}(K \cup \{y_1, y_2, \dots\}, \{F_n\}) \\
&= \max \left\{ \sup_{i \in \mathbb{N}} h_{\mathrm{top}}(K_i \cup \{y_i\}, \{F_n\}),\ h_{\mathrm{top}}(\{x_0, y_1, y_2, \dots\}, \{F_n\}) \right\} \\
&= 0,
\end{align*}
contradicting $h_{\mathrm{top}}(K, \{F_n\}) > 0$.

\textbf{Case 2}: $z \ne x_0$. 
Then there exists $I \in \mathbb{N}$ such that $z \in X_I$. Let $r = \frac{1}{2} d(z, x_0) > 0$. Since $X_i \to \{x_0\}$ in the Hausdorff metric, there exists $I_1 \in \mathbb{N}$ such that for all $i > I_1$,
\[
d(z, X_i) \ge d(z, x_0) - d_H(X_i, \{x_0\}) > r.
\]
As $z$ is the unique limit point of $K$, the set $\left( \bigcup_{i > I_1} X_i \cup \{x_0\} \right) \cap K$ is finite. Hence,
\[
h_{\mathrm{top}}(K, \{F_n\}) = h_{\mathrm{top}}\left( K \cap \left( \textstyle\bigcup_{1 \le i \le I_1} X_i \right), \{F_n\} \right) \le \max_{1 \le i \le I_1} h_{\mathrm{top}}(X_i, \{F_n\}) < h,
\]
again a contradiction.

Therefore, $\{i \in \mathbb{N} : x_0 \notin \Gamma_{\varepsilon_i}(x_i)\}$ must be finite. Without loss of generality, assume $x_0 \in \Gamma_{\varepsilon_i}(x_i)$ for all $i \ge 1$. Since $\Gamma_{\varepsilon_i}(x_i) \subseteq \Gamma_{2\varepsilon_i}(x_0)$, we have
\[
\lim_{\varepsilon \to 0^+} h_{\mathrm{top}}(\Gamma_{\varepsilon, \{F_n\}}(x_0), \{F_n\}) = \lim_{\varepsilon \to 0^+} \sup_{x \in X} h_{\mathrm{top}}(\Gamma_{\varepsilon}(x), \{F_n\}) = h.
\]

Now consider the orbit $\mathrm{orb}(x_0, G) = \{g x_0 : g \in G\}$. If it is infinite, choose $y = \lim_{j \to \infty} y_j$ where $y_j = g_{j} x_0$ with $g_{j} \in G$ and $d(y, y_j) < \frac{1}{j}$. For each $j \in \mathbb{N}$,
\[
\lim_{\varepsilon \to 0^+} h_{\mathrm{top}}(\Gamma_{\varepsilon}(y_j), \{F_n\}) = h.
\]
Select $\eta_j \in (0, d(y, y_j))$ such that
\[
h_{\mathrm{top}}(\Gamma_{\eta_j}(y_j), \{F_n\}) > \min\left\{ h(1 - \tfrac{1}{j}),\ j \right\}, \quad \forall j \in \mathbb{N}.
\]
Then $\lim_{j \to \infty} y_j = y$, $\lim_{j \to \infty} h_{\mathrm{top}}(\Gamma_{\eta_j}(y_j), \{F_n\}) = h$, and $y \notin \Gamma_{\eta_j}(y_j)$ for all $j$. As before, this leads to a contradiction.

Hence, $\mathrm{orb}(x_0, G)$ is finite. Write $\mathrm{orb}(x_0, G) = \{x_0, x_1, \dots, x_l\}$. Define
\[
Y_\varepsilon := \bigcup_{k=0}^l \mathrm{cl} \left(\Gamma_{\varepsilon}(x_k) \right).
\]
Then $Y_\varepsilon$ is a $G$-invariant closed subset, so $(Y_\varepsilon, G)$ is a subsystem of $(X,G)$ and
\[
h_{\mathrm{top}}(Y_{\varepsilon},G) \ge h_{\mathrm{top}}(\Gamma_{\varepsilon}(g_k x_0), \{F_n\}) = h .
\]
By the classical variational principle, there exist ergodic $G$-invariant probability measures $\mu_j$ on $Y_{1/j}$ such that
\[
h_{\mu_j}(Y_{1/j}, G) > \min\left\{ h(1 - \tfrac{1}{j}),\ j \right\} \quad \text{for all } j.
\]
Since $\mu_j(\Gamma_{1/j}(x_0)) > 0$, by Theorem 3.10 in \cite{Wang2026}, there exist $K_j \subseteq \Gamma_{1/j}(x_0)$ with $\mu_j(K_j) > 0$, $x_0 \notin K_j$, and
\[
h_{\mathrm{top}}(K_j, \{F_n\}) \ge h_{\mu_j}(Y_{1/j}, G) > \min\left\{ h(1 - \tfrac{1}{j}),\ j \right\}.
\]
By hereditary uniform lowerability, for each $j$ there exists an infinite countable compact subset $X_j \subseteq K_j$ with unique limit point $y_j$ in $X$ such that
\[
h_{\mathrm{top}}(X_j, \{F_n\}) = b_j := \min\left\{ h(1 - \tfrac{1}{j}),\ j \right\} < h.
\]
Then $\lim_{j \to \infty} y_j = x_0$, $\lim_{j \to \infty} h_{\mathrm{top}}(X_j, \{F_n\}) = h$, and $x_0 \notin X_j$ for all $j$. As before, this leads to a contradiction.

We conclude that
\[
\lim_{\varepsilon \to 0^+} \sup_{x \in X} h_{\mathrm{top}}(\Gamma_{\varepsilon, \{F_n\}}(x), \{F_n\})= \lim_{\varepsilon \to 0^+} \sup_{x \in X} h_{\mathrm{top}}(\Gamma_{\varepsilon}(x), \{F_n\}) = 0,
\]
so $(X,G)$ is asymptotically $h$-expansive. This completes the proof.
\end{proof}

\subsection{Hereditary uniform lowerability of asymptotically $h$-expansive $G$-systems}\label{subsection-asyhex=>hul}

Recall the key proposition proved by Dou, Zheng, and Zhou \cite[Theorem 4.2]{Dou2023ZhengZhou}.

\begin{proposition}\label{prop-Bowentype-formula}
Let $\pi : (Y,G) \to (X,G)$ be a factor map between $G$-systems, and let $\{F_n\} $ be a tempered F{\o}lner sequence of $G$.
For every non-empty compact subset $K \subseteq Y$,
\[
h_{\mathrm{top}}(\pi(K), \{F_n\}) \le h_{\mathrm{top}}(K, \{F_n\}) \le h_{\mathrm{top}}(\pi(K), \{F_n\}) + h_{\mathrm{top}}(G, Y|\pi).
\]
\end{proposition}

\begin{lemma}\label{lem-entropyproperty-principalmap}
Let $\pi : (Y,G) \to (X,G)$ be a principal factor map between $G$-systems, and let $\{F_n\}_{n \in \mathbb{N}}$ be a tempered F{\o}lner sequence of $G$. If $(X,G)$ has finite topological entropy, then 
\begin{enumerate}
\item $h_{\mathrm{top}}(G, Y|\pi) = 0$, and
\item for every non-empty compact subset $K \subseteq Y$, 
\[
h_{\mathrm{top}}(K, \{F_n\}) = h_{\mathrm{top}}(\pi(K), \{F_n\}).
\]
\end{enumerate}
\end{lemma}

\begin{proof}
Since $\pi$ is principal, Proposition~\ref{prop-VP-yan} implies that $h_{\mathrm{top}}(G, Y|\pi) = 0$. Combining this with Proposition~\ref{prop-Bowentype-formula}, we conclude that for every non-empty compact subset $K \subseteq Y$,
\[
h_{\mathrm{top}}(K, \{F_n\}) = h_{\mathrm{top}}(\pi(K), \{F_n\}).
\]
\end{proof}

We immediately obtain the following analogue of Bowen's theorem \cite[Theorem 17]{Bowen1971} for amenable group actions.

\begin{theorem}
Let $\pi : (Y,G) \to (X,G)$ be a factor map between $G$-systems. If $(X,G)$ has finite topological entropy, then
\[
h_{\mathrm{top}}(Y, G) \le h_{\mathrm{top}}(X, G) + h_{\mathrm{top}}(G, Y|\pi).
\]
\end{theorem}

\begin{proof}
Let $\{F_n\} $ be a tempered F{\o}lner sequence of $G$. Applying Proposition~\ref{prop-Bowentype-formula} with the tempered F{\o}lner sequence $\{F_n\} $ and $K = Y$, we obtain
\[
h_{\mathrm{top}}(Y, \{F_n\}) \le h_{\mathrm{top}}(\pi(Y), \{F_n\}) + h_{\mathrm{top}}(G, Y| \pi).
\]
Since $\pi(Y) = X$ and both the topological entropy of a $G$-system and the relative topological entropy with respect to $\pi$ are independent of the choice of F{\o}lner sequences, the result follows.
\end{proof}

We next prove the other direction of the equivalence in Theorem~\ref{thm-main-HUL}.

\begin{theorem}\label{thm-asyhex->hul}
Let $(X,G)$ be an asymptotically $h$-expansive $G$-system, and let $\{F_n\}$ be a tempered F{\o}lner sequence of $G$ satisfying the conditions (\ref{FS-condition}). Then $(X,G)$ is hereditarily uniformly lowerable along $\{F_n\}$, i.e., for every non-empty compact subset $K \subseteq X$ and every $0 \le h \le h_{\mathrm{top}}(K, \{F_n\})$, there exists a non-empty compact subset $K_h \subseteq K$ such that $h_{\mathrm{top}}(K_h, \{F_n\}) = h$ and $K_h$ has at most one limit point.
\end{theorem}

\begin{proof}
By Lemma~\ref{lem-equicondition-asyhexpan}, there exists a principal extension $\pi: (\overline{Y}, G) \to (X,G)$ such that $(\overline{Y}, G)$ is a quasi-symbolic system, i.e., $\overline{Y} = Y \vee \mathbf{T}$, where $Y$ is a subshift and $\mathbf{T}$ is a tiling system with zero entropy.

Clearly, $\overline{Y}$ is a principal extension of $Y$, denoted by $\pi_Y: \overline{Y} \to Y$. Since $\mathbf{T}$ has zero topological entropy, we have
\[
h_{\mathrm{top}}(X, G) \le h_{\mathrm{top}}(\overline{Y}, G) = h_{\mathrm{top}}(Y, G) < \infty.
\]

Fix a non-empty compact subset $K \subseteq X$ and $0 \le h \le h_{\mathrm{top}}(K, \{F_n\})$. Since $\pi$ and $\pi_Y$ are principal, Lemma~\ref{lem-entropyproperty-principalmap} implies that $K^Y := \pi_Y(\pi^{-1}(K)) \subseteq Y$ is a non-empty compact subset satisfying
\[
h_{\mathrm{top}}(K^Y, \{F_n\}) = h_{\mathrm{top}}(K, \{F_n\}).
\]

By Theorem~\ref{thm-expansive-hul}, for the given $0 \le h \le h_{\mathrm{top}}(K^Y, \{F_n\})$, there exists an infinite countable compact subset $K_h^Y \subseteq K^Y$ such that
\[
h_{\mathrm{top}}(K_h^Y, \{F_n\}) = h \quad \text{and} \quad K_h^Y \text{ has at most one limit point in } Y.
\]
Write $K_h^Y = \{y_0, y_1, y_2, \dots\}$ with $y_i \to y_0$ as $i \to \infty$.

We now utilize the specific construction of Downarowicz and Zhang in the proof of Theorem 5.2 in \cite{Downarowicz2023Zhang}.
In their proof, $\overline{Y}$ is constructed from a principal zero-dimensional extension $\overline{X}$ of $X$, denoted by $\pi_X: \overline{X} \to X$, where $\overline{X}$ has $\mathbf{T}$ as a topological factor. The space $\overline{Y}$ is obtained as an intersection of spaces $\overline{Y}_k \subseteq Y_k \vee \mathbf{T}$, each factoring via a map $\overline{\pi}_k$ onto $\overline{X}_{[1,k]}$, with these factor maps being consistent. The factor map $\overline{\pi}: \overline{Y} \to \overline{X}$ is then defined as the limit of the block codes $\overline{\pi}_k$.

For each $i \in \mathbb{N} \cup \{0\}$ and $k \in \mathbb{N}$, define
\begin{align*}
\mathbf{T}_k(y_i) &= \left\{ \mathcal{T}_k \in \mathbb{T}_k : \pi_Y(y_i, \mathcal{T}) = y_0 \text{ for some } \mathcal{T} = (\mathcal{T}_k)_{k \in \mathbb{N}} \in \mathbf{T} \right\} \\
&= \left\{ \mathcal{T}_k \in \mathbb{T}_k : (y_i, \mathcal{T}) \in \overline{Y} \text{ for some } \mathcal{T} = (\mathcal{T}_k)_{k \in \mathbb{N}} \in \mathbf{T} \right\},
\end{align*}
where $ \mathbf{T} = \stackrel{\longleftarrow}{\lim\limits_{k}} \mathbb{T}_k $. 

By the construction of $\overline{Y}_k$, there exists $N(1) \in \mathbb{N}$ such that for all $i \ge N(1)$, we have $\mathbf{T}_1(y_i) = \mathbf{T}_1(y_0)$. For each $k \ge 2$, there exists $N(k) \in \mathbb{N}$ with $N(k) > N(k-1)$ such that for all $i \ge N(k)$, we have $\mathbf{T}_k(y_i) = \mathbf{T}_k(y_0)$.

For each $(y_0, \mathcal{T}) \in \pi_Y^{-1}(y_0)$ and $k \in \mathbb{N}$, note that
\[
(\overline{\pi}_{k+1})^{-1}(\overline{\pi}(y_0, \mathcal{T})_{[1,k+1]}) \subseteq (\overline{\pi}_k)^{-1}(\overline{\pi}(y_0, \mathcal{T})_{[1,k]}),
\]
where $\overline{\pi}(y_0, \mathcal{T})_{[1,k]}$ denotes the projection onto the first $k$ layers of $\overline{X}$.

We may therefore select $\overline{y}_0 = (y_0, \mathcal{T}) \in \pi_Y^{-1}(y_0)$ with $\mathcal{T} = (\mathcal{T}_k)_{k \in \mathbb{N}}$ and $\mathcal{T}_k \in \mathbf{T}_k(y_0)$ for each $k \in \mathbb{N}$.
Since $y_i \to y_0$ as $i \to \infty$, for each $k \in \mathbb{N}$ and $N(k) \le i < N(k+1)$, we can choose $\overline{y}_i = (y_i, \mathcal{T}^i) \in \pi_Y^{-1}(y_i)$ such that 
\[ 
\mathcal{T}^{i}_l = \mathcal{T}_l \quad \text{for all } 1 \le l \le k .\]
For $i < N(1)$, we choose $\overline{y}_i = (y_i, \mathcal{T}^i) \in \overline{Y}$ arbitrarily. Then $\overline{y}_i \to \overline{y}_0$ as $i \to \infty$.

Define $K_h^{\overline{Y}} = \{\overline{y}_i : i \in \mathbb{N}\}$. Then $\pi_Y(K_h^{\overline{Y}}) = K_h^Y$, and $K_h^{\overline{Y}}$ is an infinite countable compact set with a unique limit point in $\overline{Y}$.
By Lemma~\ref{lem-entropyproperty-principalmap},
\[
h_{\mathrm{top}}(K_h^{\overline{Y}}, \{F_n\}) = h_{\mathrm{top}}(K_h^Y, \{F_n\}) = h.
\]
Finally, set $K_h = \pi(K_h^{\overline{Y}})$. Then $K_h$ is a countable compact subset of $X$ with at most one limit point, and
\[
h_{\mathrm{top}}(K_h, \{F_n\}) = h_{\mathrm{top}}(K_h^{\overline{Y}}, \{F_n\}) = h,
\]
which shows that $(X,G)$ is hereditarily uniformly lowerable along $\{F_n\}$.
\end{proof}

\begin{proof}[Proof of Theorem~\ref{thm-main-HUL}]
	The two implications established in Theorem~\ref{thm-hul->asyhex} and Theorem~\ref{thm-asyhex->hul} complete the proof.
\end{proof}

\subsection{Lowerability under principal extensions}\label{subsection-preservation}

\begin{proposition}
Let $\pi :(Y,G) \to (X,G)$ be a principal factor map between $G$-systems, and let $\{F_n\}$ be a tempered F{\o}lner sequence of $G$ satisfying the
conditions (\ref{FS-condition}).
Then $(Y,G)$ is lowerable (respectively, hereditarily lowerable, hereditarily uniformly lowerable) if and only if $(X,G)$ has the same property.
\end{proposition}

\begin{proof}
The preservation of the lowerable and hereditarily lowerable  under a principal extension can be obtained directly by Lemma \ref{lem-entropyproperty-principalmap}. It remains to verify the case for hereditarily uniformly lowerable.

First, suppose $(Y,G)$ is hereditarily uniformly lowerable along $\{F_n\}$. For any non-empty compact subset $K \subseteq X$ and $0 \le h \le h_{\mathrm{top}}(K, \{F_n\})$, let $K_Y = \pi^{-1}(K)$. By Lemma~\ref{lem-entropyproperty-principalmap}, we have
\[
h_{\mathrm{top}}(K_Y, \{F_n\}) = h_{\mathrm{top}}(K, \{F_n\}).
\]
Since $(Y,G)$ is hereditarily uniformly lowerable, there exists a non-empty compact subset $K_{Y,h} \subseteq K_Y$ such that
\[
h_{\mathrm{top}}(K_{Y,h}, \{F_n\}) = h \quad \text{and} \quad K_{Y,h} \text{ has at most one limit point in } Y.
\]
Set $K_h = \pi(K_{Y,h})$. Then
\[
h_{\mathrm{top}}(K_h, \{F_n\}) = h \quad \text{and} \quad K_h \text{ has at most one limit point in } X.
\]
This shows that $(X,G)$ is hereditarily uniformly lowerable along $\{F_n\}$.

Conversely, suppose $(X,G)$ is hereditarily uniformly lowerable along $\{F_n\}$. By Theorem~\ref{thm-main-HUL}, $(X,G)$ is asymptotically $h$-expansive, i.e., $h_{\mathrm{top}}^*(X, G) = 0$. Combining this with Lemma~\ref{lem-relative&conditional} and Proposition~\ref{prop-VP-yan}, we obtain
\begin{align*}
h_{\mathrm{top}}(G,Y|\pi) = h_{\mathrm{top}}^*(G,Y|\pi)=0.
\end{align*}
Then by Lemma \ref{lem-condition-asyhex}, 
\begin{align*}
h_{\mathrm{top}}^*(Y,G) = 0,
\end{align*}
which means $(Y,G)$ is asymptotically $h$-expansive. Applying Theorem~\ref{thm-main-HUL} again, we conclude that $(Y,G)$ is hereditarily uniformly lowerable along $\{F_n\}$.
This completes the proof.
\end{proof}

In fact, the preceding proposition implies the following result, since the conditional topological entropy is independent of the choice of F{\o}lner sequences of $G$.

\begin{proposition}
Let $\pi : (Y,G) \to (X,G)$ be a principal factor map between $G$-systems. Then $(Y,G)$ is asymptotically $h$-expansive if and only if $(X,G)$ is asymptotically $h$-expansive.
\end{proposition}


$\mathbf{Acknowledgements}$.The authors wish to express their gratitude to Professor Dou Dou for his insightful comments and constructive discussions.

\bibliographystyle{amsplain}

\end{document}